\documentclass[manuscript,screen]{acmart}

\usepackage{algorithm}
\usepackage{algorithmic}
\usepackage{listings}
\usepackage{tikz}
\usetikzlibrary{calc,shapes}
\usepackage{array}
\usepackage{pgfplots}
\usepackage{cleveref}
\usepackage{pythonhighlight}
\usepackage{usebib}
\usepackage{mathtools}
\usepackage{ifthen}

\newcommand{\anonymous}{NO}

\copyrightyear{2026}
\acmYear{2026}
\acmDOI{XXXXXXX.XXXXXXX}

\newtheorem{definition}{Definition}

\newtheorem{assumption}{Assumption}

\crefname{figure}{figure}{figures}
\Crefname{figure}{Figure}{Figures}
\crefname{assumption}{assumption}{assumptions}
\Crefname{assumption}{Assumption}{Assumptions}
\crefname{theorem}{example}{example}
\Crefname{theorem}{Example}{Examples}

\newcommand{\refel}{R}
\newcommand{\polyspace}{\mathcal{V}}
\newcommand{\dualbasis}{\mathcal{L}}
\newcommand{\dual}[1]{#1^*}

\newcommand{\transpose}{^{\mathsf{T}}}
\newcommand{\invtranspose}{^{-\mathsf{T}}}

\newcommand{\mat}[1]{\boldsymbol{\mathrm{#1}}}
\renewcommand{\vec}[1]{\boldsymbol{#1}}

\newcommand{\TODO}[1][]{{\color{red}TODO\ifthenelse{\equal{#1}{}}{}{: #1}}}

\newcommand{\rtriangle}{\hspace{1pt}\begin{tikzpicture}[x=1.5mm,y=1.5mm]\draw (0,0) -- (1,0) -- (0,1) -- cycle;\end{tikzpicture}}
\newcommand{\interval}{\hspace{1pt}\raisebox{0.8mm}{\begin{tikzpicture}[x=1.5mm,y=1.5mm]\draw (0,0) -- (1,0) (0,0.3) -- (0,-0.3) (1,0.3) -- (1,-0.3);\end{tikzpicture}}}
\newcommand{\quadrilateral}{\hspace{1pt}\begin{tikzpicture}[x=1.5mm,y=1.5mm]\draw (0,0) -- (1,0) -- (1, 1) -- (0,1) -- cycle;\end{tikzpicture}}
\newcommand{\twodimentity}{\hspace{1pt}\begin{tikzpicture}[x=1.5mm,y=1.5mm]\draw (0,1) -- (0,0) -- (1,0);\foreach \x in {0.75,0.5,0.25}\fill (\x+0.1,1-\x+0.1) circle (0.3pt);\end{tikzpicture}}
\newcommand{\reflection}{\textup{ref}}
\newcommand{\rotation}{\textup{rot}}
\newcommand{\scaleddown}[1]{\scalebox{0.8}{#1}}

\newcommand{\prismname}{triangular prism}
\newcommand{\pyramidname}{square-based pyramid}

\newcommand{\modifiedfunctional}{\bar{l}}
\newcommand{\modifiedbasis}{\bar{\phi}}

\usepackage{xspace,color}

\begin{document}

\title{Algorithm XXXX: Computation of finite element degree-of-freedom transformation matrices}


\ifthenelse{\equal{\anonymous}{YES}}{
\author{Author Name}
\orcid{ANON-YMOU-SORC-IDxx}
\email{anon@example.com}
\affiliation{%
  \institution{University}
  \country{Country}
}
\author{Author Name}
\orcid{ANON-YMOU-SORC-IDxx}
\email{anon@example.com}
\affiliation{%
  \institution{University}
  \country{Country}
}
}{ 
\author{Matthew W.~Scroggs}
\orcid{0000-0002-4658-2443}
\email{matthew.scroggs.14@ucl.ac.uk}
\affiliation{%
  \institution{Advanced Research Computing Centre, University College London}
  \streetaddress{Gower Street}
  \city{London}
  \country{United Kingdom}
  \postcode{WC1E 6BT}
}
\author{Garth N.~Wells}
\email{gnw20@cam.ac.uk}
\orcid{0000-0001-5291-7951}
\affiliation{%
  \institution{Department of Engineering, University of Cambridge}
  \streetaddress{Trumpington Street}
  \city{Cambridge}
  \country{United Kingdom}
  \postcode{CB2 1PZ}
}
}

\begin{abstract}
The arithmetic intensity of algorithms for computing finite element 
operators increases with increasing polynomial degree. This has made 
high degree methods particularly attractive on modern CPU and GPU 
architectures, since performance at low degree is limited by memory 
bandwidth and only a small fraction of the floating point capacity of 
the processor is used. Higher degree methods can exploit a 
significantly greater fraction of the available compute power. While 
stable methods for computing high-degree finite element bases are 
well-established, there is no universal and automated algorithm for 
the efficient construction of the degree-of-freedom map for arbitrary 
degree elements. We introduce a new algorithm that can be used in 
computing degree-of-freedom maps for arbitrary Ciarlet-type finite 
elements using only the element's definition and properties of the 
reference cell, without requiring a specific implementation for each 
element. This method is implemented in the library Basix, a component 
of the FEniCSx libraries. As well as allowing vast simplifications of 
parts of a codebase, the algorithm has allowed us to support 
user-defined custom elements that a user can create at runtime without 
requiring the user to input any information about transformations 
required to construct a degree-of-freedom map.
\end{abstract}



\maketitle

\section{Introduction}

In finite element libraries it is usual for global finite element
vectors or matrices to be computed by evaluating cell-wise contributions
and combining these to form a global vector or matrix. The scattering of
cell-wise contributions to the global vector/matrix must preserve the
required continuity of finite element functions between cells. The
local-to-global map that ensures this continuity is often referred to as
the \emph{degree-of-freedom map}.

\begin{figure}
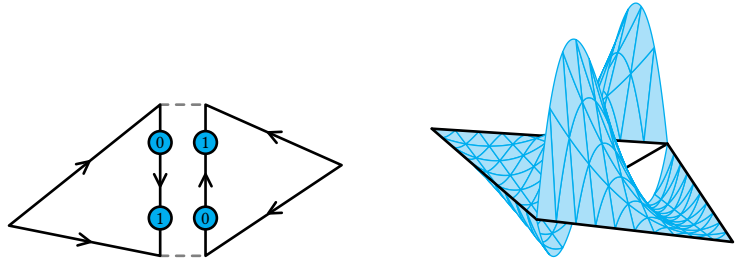



\caption{In an unstructured mesh, two neighbouring cells will not
necessarily agree on the orientation of their sub-entities: the two
triangles shown in the left diagram orient their shared edge differently
(the triangles are shown here with a small gap between them and a
duplicated edge to allow for clear labelling; in the actual mesh, the
edge will be shared). This causes to the two cells to place the 0th and
1st DOFs associated with this edge in a degree 3 Lagrange space at
different points, leading to the incorrect basis function shown in the
right diagram.}
\label{fig:mismatched_edge}
\end{figure}

Degrees-of-freedom (DOFs) of an element can be associated with cell
(sub-)entities, i.e.~vertices, edges, faces or the cell volume. When
using higher-degree finite element spaces, there can be multiple DOFs
associated with sub-entities that are shared by more than one cell
(e.g., in a degree 3 Lagrange space on a triangle or quadrilateral,
there are two DOFs associated with each edge, and edges can be shared by
two cells). To ensure the required continuity between cells,
neighbouring cells must agree on the orientation of shared sub-entities.
Failure to do this can lead to a mismatch in the arrangement of the DOFs
on shared sub-entities and incorrect combinations of values being
inserted into the global matrix (see \cref{fig:mismatched_edge}).
Agreement on a common orientation is not limited to elements with more
than one DOF associated with a cell entity; it is also required for
elements with DOFs that are defined in terms of orientation,
e.g.~$H(\operatorname{div})$- and $H(\operatorname{curl})$-conforming
finite elements. An approach to agreement on common entity orientations
is through global mesh orderings that ensure that guarantee this
property (see \cite{2022-dofs} for an extensive discussion of the
published approaches). However, not all meshes of hexahedral cells can
be suitably ordered \cite{agelek:2017} and meshes of mixed cells types
pose particular challenges.

In \cite{2022-dofs}, we presented a method of permutations and
transformations that can be applied to the local contributions from each
cell to correct for any orientation differences in an unstructured mesh.
The method was implemented in the library Basix \cite{basix}, the finite
element definition and tabulation library used by FEniCSx
\cite{fenicsx}. Our initial implementation involved bespoke code for
each element type to generate a set of `base transformation matrices'.
In this paper, we present a new algorithm that computes the base
transformations from the definition of the element and properties of the
reference cell that it is defined on. This algorithm is general and does
not require any specific implementation for each element. As well as
allowing us to greatly simplify the implementation in Basix, this new
algorithm allows us to support user-defined custom elements
\cite[section~5]{fenicsx} without requiring the user to provide
complicated information on the DOF transformations.

A high-performance C++ implementation of the developed algorithm was
added to Basix in 2022 \cite{basix-pr}. An alternative method for
algorithmically computing matrices for this purpose was independently
proposed by T.~Collin in 2025 \cite[section 6.4]{collin-thesis}.
In this paper's software component \cite{supplement}, we provide 
the Basix implementation of the algorithm alongside a simpler
demonstration implementation in Python using the symbolic finite element prototyping
library Symfem \cite{symfem}. The source code of both of these
implementations is available under the MIT open source license.

The remainder of the paper is laid out as follows. In
\cref{sec:finite-element}, we present the Ciarlet definition of a finite
element, on which our algorithm is built. In
\cref{sec:dof-transformations}, we outline our method of DOF
transformations and describe how they can be represented using a small
set of `base transformation matrices'. In \cref{sec:algorithm}, we
describe the new algorithm for computing the base transformation
matrices, and  in \cref{sec:in-place} we present a method of in-place
matrix multiplication that can be used to efficiently apply the base
transformation matrices. We finish with some concluding remarks in
\cref{sec:outro}.

\section{Defining a finite element}
\label{sec:finite-element}

In general, finite elements can be defined as
follows~\cite{Ciarlet:1978}.
\begin{definition}[Ciarlet finite element]\label{ciarlet-def}
A finite element is defined by the triple
$(\refel,\polyspace,\dualbasis)$, where
\begin{itemize}
  \item $\refel\subset\mathbb{R}^{d_\refel}$ is the reference cell,
  usually a polygon or polyhedron;
  \item $\polyspace$ is a finite dimensional space on
  $\refel$ of dimension~$n$, usually a space of polynomials;
  \item $\dualbasis \coloneq\{\hat{l}_0,\dots,\hat{l}_{n-1}\}$ is a basis of
  the dual space $\dual{\polyspace} : =
  \{f:\polyspace\to\mathbb{R}\,|\,f\text{ is linear}\}$. Each functional
  $\hat{l}_i$ is associated with a sub-entity of the reference cell
  $\refel$.
\end{itemize}
The reference basis functions $\{\hat{\phi}_0, \dots,
\hat{\phi}_{n-1}\}$ of the space $\polyspace$ are defined by
\begin{equation*}
  \hat{l}_i(\hat{\phi}_j)= \delta_{ij}\coloneq\begin{cases}1 & i=j,
    \\
    0 & i\not=j.\end{cases}
\end{equation*}
If the functional $\hat{l}_i$ is associated with the sub-entity $E$,
then we may also say that the basis function $\hat{\phi}_i$ `is
associated with~$E$'.
\end{definition}
The value $d_\refel$ is the topological dimension of the cell; this may
differ from the geometric dimension $d$ of the cell if, for example, a
mesh of polygon cells is embedded in $\mathbb{R}^3$. Further, a map is
required that maps basis functions on the reference cell to functions on
a physical cell \cite{rognes:2009,mapping2,mapping3}. We refer to the
basis function map from the reference cell to a physical cell as the
\emph{push-forward}. For a geometry map $g : \refel \to \mathbb{R}^d$,
we denote the push-forward map by $\mathcal{F}_g : \polyspace \to
\mathcal{W}$, where $\mathcal{W}$ is a function space on the physical
cell $g(\refel)$. The map $\mathcal{F}_g$ for an element is chosen such
that required properties of the basis functions are preserved on the
physical cell.
The functionals $\hat{l}_i \in \dualbasis$ are the (local)
degrees-of-freedom (DOFs) of the finite element. Note that we enumerate
functionals and basis functions from~0. When a finite element function
space is defined on a mesh, we associate a global DOF index with each
local DOF on each cell. To ensure that the mapped space has the required
continuity properties, any local DOF that is associated with a
sub-entity that is shared by multiple cells must be assigned the same
global DOF number as the corresponding DOF on the neighbouring cell(s).

\label{entity-by-entity}
In Basix, the DOFs on each element are numbered entity-by-entity; the
DOFs for each sub-entity are contiguously numbered, with those for
sub-entities of lower dimension appearing first; for sub-entities of the
same dimension, those associated with the sub-entity of the lower index
have the lower indices. For simplicity, in this paper we assume that
this entity-by-entity numbering is used, although the methods presented
can be adapted to other ordering conventions.

We make the following assumptions about the functionals associated with
each sub-entity, which are true of all widely used finite element
spaces, and are key to implementations being able to enforce continuity
between neighbouring cells.
\begin{assumption}\label{assume-me}
Each sub-entity of the same type (e.g.~each sub-entity that is a
triangle) has an equivalent set of DOF functionals associated with it.
\end{assumption}
\begin{assumption}\label{assume-me2}
Each functional $\hat{l}_i \in \dualbasis$ depends only on the values of
functions restricted to the sub-entity with which $\hat{l}_i$ is
associated.
\end{assumption}
\begin{assumption}\label{assume-me3}
Let $E$ be a sub-entity of $\refel$, let $g:E\to E$ be an affine
bijection, and $\mathcal{F}_g: \left.\polyspace \right|_E \to \left.
\mathcal{W} \right|_E$ be the push forward corresponding to $g$. If
$\hat{\phi}_a,\dots,\hat{\phi}_b$ are the basis functions associated
with the sub-entity $E$, then
$$
\operatorname{span}\left(\left\{\mathcal{F}_g\left(\left.\hat{\phi}_a\right|_E\right),
    \dots, \mathcal{F}_g\left(\left.\hat{\phi}_b\right|_E\right)\right\}\right)
=
\operatorname{span}\left(\left\{\left.\hat{\phi}_a\right|_E,
    \dots, \left.\hat{\phi}_b\right|_E\right\}\right),
$$
where $|_E$ denotes the restriction of a function or function space to
the sub-entity~$E$.
\end{assumption}
\Cref{assume-me} ensures that two neighbouring cells will have equivalent
DOFs on their shared sub-entities. \Cref{assume-me2} ensures that two
neighbouring cells only need to `communicate' using function values on
shared sub-entities. By \cite[lemma~20.6, exercise~20.1]{ern-guermond},
we see that any map that corresponds to a rotation or reflection of a
sub-entity is affine and a bijection, hence \cref{assume-me3} ensures that
if you have two neighbouring cells and you rotate or reflect one of
them, the basis functions associated a sub-entity shared by both cells
will span the same space restricted to that entity.

\subsection{Defining functionals via quadrature}
\label{sec:basix-functionals}

In Basix, the implementation of elements is based on the Ciarlet
definition (\cref{ciarlet-def}). The simplest type of functional is a
point evaluation functional, which is applied to a function $f$ by
evaluating $f$ at a given point on the reference cell. Point evaulation
functionals are used when defining Lagrange elements. For other finite
elements---including Raviart--Thomas \cite{rt}, N\'ed\'elec
\cite{nedelec,nedelec2}, Brezzi--Douglas--Marini \cite{bdm} and
serendipity \cite{serendipity} elements, and more
\cite{commonandunusual,defelement,defelement-paper}---the functionals in
$\dualbasis$ include integral moments. These integral moment functionals
are evaluated by multiplying a function $f$ by a given weight function
and integrating over a sub-entity of the cell. When implementing finite
elements, these integrals can be computed using numerical quadrature.

In Basix, a functional associated with a given sub-entity $E$ is defined
by (i) a set of points $\{\vec{p}_0, \dots, \vec{p}_{n_p-1}\} \in
\mathbb{R}^{d_\refel}$, where $n_p$ is the number of points, and (ii) a
weight tensor $\mat{M}=(m_{ijk}) \in \mathbb{R}^{n_E \times s \times
n_p}$ (where $n_E$ is the number of DOFs associated with the sub-entity
$E$ and $s$ is the value size of the element). The functionals
$\hat{l}_a, \dots, \hat{l}_b$ associated with $E$ are evaluated using
$$
  \hat{l}_{a+i}(f)
  = \sum_{j=0}^{s-1} \sum_{k=0}^{n_p-1} m_{ijk} \left[f(\vec{p}_k )\right]_j,
$$
where $\left[f\right]_j$ denotes the $j$th component of $f$ if $f$ is a
vector-valued function; if $f$ is a scalar-valued function, then $\left[
f \right]_0 \coloneq f$. Functionals of various types can be implemented using
this representation: for point evaluation DOFs, we can use a single
point and a single weight of 1, and for integral moments we can use the
quadrature points and weights. As the functionals associated with each
sub-entity of the cell are usually of the same type, we include one set
of points $\{\vec{p}_0, \dots\}$ per sub-entity and use these for each
functional associated with that sub-entity. Where functionals associated
with a sub-entity have a mixture of types, points that define both can
be included, with the weights corresponding to points not relevant to a
given functional set to~0.

\subsection{Examples}

Before we describe our algorithm, we consider the definitions of some
common elements. These elements will later be used in the examples of
the computation of DOF transformations. The diagrams in this section are
taken from DefElement~\cite{defelement,defelement-paper}. In this
section and throughout this paper, we use the same definitions of
reference cells as used in Basix (which are the same as those used by
Symfem and DefElement, although the numbering of sub-entities for simplex
cells is done differently to Basix), but we note that these are arbitrary and the
methods presented could be adapted to any choice of reference cell.

\begin{figure}
\begin{tikzpicture}[x=1.5cm,y=1.5cm,line cap=round,line join=round,line width=1pt]
\draw[-stealth,black] (25.0,25.0) -- (25.375,25.0);
\draw[-stealth,black] (25.0,25.0) -- (25.0,25.375);
\node[black,anchor=west] at (25.35,25.0) {\footnotesize$x$};\node[black,anchor=south] at (25.0,25.35) {\footnotesize$y$};\draw[black] (26.125,25.0) -- (27.25,25.0);
\draw[black] (26.125,25.0) -- (26.125,26.125);
\draw[black] (27.25,25.0) -- (27.25,26.125);
\draw[black] (26.125,26.125) -- (27.25,26.125);
\draw[orange,fill=white] (26.125,25.0) circle (4.0pt);
\node[black,anchor=center] at (26.125,25.0) {\footnotesize 0};\draw[orange,fill=white] (27.25,25.0) circle (4.0pt);
\node[black,anchor=center] at (27.25,25.0) {\footnotesize 1};\draw[orange,fill=white] (26.125,26.125) circle (4.0pt);
\node[black,anchor=center] at (26.125,26.125) {\footnotesize 2};\draw[orange,fill=white] (27.25,26.125) circle (4.0pt);
\node[black,anchor=center] at (27.25,26.125) {\footnotesize 3};\draw[black] (27.625,25.0) -- (28.75,25.0);
\draw[cyan,fill=white] (28.1875,25.0) circle (4.0pt);
\node[black,anchor=center] at (28.1875,25.0) {\footnotesize 0};\draw[black] (27.625,25.0) -- (27.625,26.125);
\draw[cyan,fill=white] (27.625,25.5625) circle (4.0pt);
\node[black,anchor=center] at (27.625,25.5625) {\footnotesize 1};\draw[black] (28.75,25.0) -- (28.75,26.125);
\draw[cyan,fill=white] (28.75,25.5625) circle (4.0pt);
\node[black,anchor=center] at (28.75,25.5625) {\footnotesize 2};\draw[black] (27.625,26.125) -- (28.75,26.125);
\draw[cyan,fill=white] (28.1875,26.125) circle (4.0pt);
\node[black,anchor=center] at (28.1875,26.125) {\footnotesize 3};\draw[green,fill=white] (29.6875,25.5625) circle (4.0pt);
\node[black,anchor=center] at (29.6875,25.5625) {\footnotesize 0};\draw[black] (29.125,25.0) -- (30.25,25.0);
\draw[black] (29.125,25.0) -- (29.125,26.125);
\draw[black] (30.25,25.0) -- (30.25,26.125);
\draw[black] (29.125,26.125) -- (30.25,26.125);
\end{tikzpicture}
\caption{The numbering of the sub-entities of a quadrilateral reference cell.}
\label{fig:quad-reference}
\end{figure}
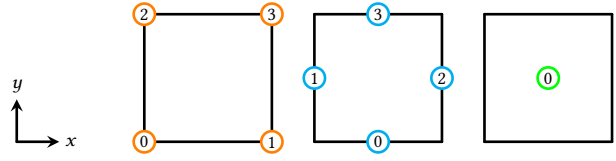

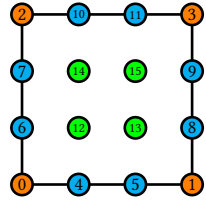
\begin{figure}
\begin{tikzpicture}[x=1cm,y=1cm,line cap=round,line join=round,line width=1pt]
\fill[white,opacity=0.5] (25.0,25.0) -- (27.25,25.0) -- (27.25,27.25) -- (25.0,27.25) -- cycle;\draw[black] (25.0,25.0) -- (27.25,25.0);
\draw[black] (25.0,25.0) -- (25.0,27.25);
\draw[black] (27.25,25.0) -- (27.25,27.25);
\draw[black] (25.0,27.25) -- (27.25,27.25);
\draw[black,fill=green] (25.75,25.75) circle (4.0pt);
\node[black,anchor=center] at (25.75,25.75) {\tiny 12};\draw[black,fill=green] (26.5,25.75) circle (4.0pt);
\node[black,anchor=center] at (26.5,25.75) {\tiny 13};\draw[black,fill=green] (25.75,26.5) circle (4.0pt);
\node[black,anchor=center] at (25.75,26.5) {\tiny 14};\draw[black,fill=green] (26.5,26.5) circle (4.0pt);
\node[black,anchor=center] at (26.5,26.5) {\tiny 15};\draw[black,fill=cyan] (25.75,25.0) circle (4.0pt);
\node[black,anchor=center] at (25.75,25.0) {\footnotesize 4};\draw[black,fill=cyan] (26.5,25.0) circle (4.0pt);
\node[black,anchor=center] at (26.5,25.0) {\footnotesize 5};\draw[black,fill=cyan] (25.0,25.75) circle (4.0pt);
\node[black,anchor=center] at (25.0,25.75) {\footnotesize 6};\draw[black,fill=cyan] (25.0,26.5) circle (4.0pt);
\node[black,anchor=center] at (25.0,26.5) {\footnotesize 7};\draw[black,fill=cyan] (27.25,25.75) circle (4.0pt);
\node[black,anchor=center] at (27.25,25.75) {\footnotesize 8};\draw[black,fill=cyan] (27.25,26.5) circle (4.0pt);
\node[black,anchor=center] at (27.25,26.5) {\footnotesize 9};\draw[black,fill=cyan] (25.75,27.25) circle (4.0pt);
\node[black,anchor=center] at (25.75,27.25) {\tiny 10};\draw[black,fill=cyan] (26.5,27.25) circle (4.0pt);
\node[black,anchor=center] at (26.5,27.25) {\tiny 11};\draw[black,fill=orange] (25.0,25.0) circle (4.0pt);
\node[black,anchor=center] at (25.0,25.0) {\footnotesize 0};\draw[black,fill=orange] (27.25,25.0) circle (4.0pt);
\node[black,anchor=center] at (27.25,25.0) {\footnotesize 1};\draw[black,fill=orange] (25.0,27.25) circle (4.0pt);
\node[black,anchor=center] at (25.0,27.25) {\footnotesize 2};\draw[black,fill=orange] (27.25,27.25) circle (4.0pt);
\node[black,anchor=center] at (27.25,27.25) {\footnotesize 3};\end{tikzpicture}
\caption{The DOFs of a degree 3 Lagrange element on a quadrilateral.}
\label{fig:quad-dofs}
\end{figure}

\begin{example}[Lagrange degree 3 on a quadrilateral with equally-spaced evaluation points]
\label{example-def-lagrange}
A degree 3 Lagrange element on a quadrilateral cell is defined by
\begin{itemize}
    \item $\refel = [0,1]^2$ (where we number the
    sub-entities as shown in \cref{fig:quad-reference}),
    \item $\polyspace = \operatorname{span} \left\{x^iy^j\, \middle|\,
    i, j \in \{0, 1, 2, 3\}\right\}$,
    \item $\dualbasis = \{\hat{l}_0, \dots, \hat{l}_{15}\}$,
where $\hat{l}_0$ to $\hat{l}_3$ are point evaluations at the vertices
of the cell (with each functional associated with its vertex);
$\hat{l}_4$ and $\hat{l}_5$ are point evaluations on edge 0 of the cell
(at points $1/3$ and $2/3$ of the way along the edge); $\hat{l}_6$ to
$\hat{l}_{11}$ are the corresponding point evaluations for edges 1 to 3;
and $\hat{l}_{12}$ to $\hat{l}_{15}$ are point evaluations on the
interior of the cell at the points $(1/3,1/3)$, $(2/3, 1/3)$, $(1/3,
2/3)$ and $(2/3, 2/3)$.
\end{itemize}
In Basix, the functionals for a sub-entity $E$ are implemented with the
relevant points and a tensor $\mat{M} = (m_{i 0 k} )\in \mathbb{R}^{n_E
\times 1 \times n_E}$, where $n_E$ is the number of points and
$$
 m_{i 0 k} =\begin{cases}1&i=k,\\0&i\not=k. \end{cases}
$$
A representation of the DOFs of this element is shown in
\cref{fig:quad-dofs}.

The push-forward for this element is the identity map. If $g:\refel \to
\mathbb{R}^d$ maps points on the reference cell to a physical cell, then
the identity map $\mathcal{F}_g^\textup{id}$ is defined by
$$
\mathcal{F}_g^\textup{id}(\hat{\psi}) = \hat{\psi}\circ g^{-1},
$$
where $\hat{\psi}$ is a function defined on the reference cell that the
push-forward is being applied to.
\end{example}

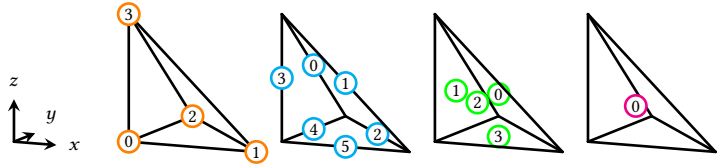
\begin{figure}
\begin{tikzpicture}[x=1.5cm,y=1.5cm,line cap=round,line join=round,line width=1pt]
\draw[-stealth,black] (25.0,25.09) -- (25.375,25.06);
\draw[-stealth,black] (25.0,25.09) -- (25.1875,25.165);
\draw[-stealth,black] (25.0,25.09) -- (25.0,25.465);
\node[black,anchor=west] at (25.405,25.0576) {\footnotesize$x$};\node[black,anchor=south west] at (25.1875,25.165) {\footnotesize$y$};\node[black,anchor=south] at (25.0,25.495) {\footnotesize$z$};\draw[black] (26.575,25.315) -- (26.0125,26.215);
\draw[black] (27.1375,25.0) -- (26.575,25.315);
\draw[black] (26.0125,25.09) -- (26.575,25.315);
\draw[orange,fill=white] (26.575,25.315) circle (4.0pt);
\node[black,anchor=center] at (26.575,25.315) {\footnotesize 2};\draw[black] (27.1375,25.0) -- (26.0125,26.215);
\draw[black] (26.0125,25.09) -- (26.0125,26.215);
\draw[black] (26.0125,25.09) -- (27.1375,25.0);
\draw[orange,fill=white] (26.0125,25.09) circle (4.0pt);
\node[black,anchor=center] at (26.0125,25.09) {\footnotesize 0};\draw[orange,fill=white] (27.1375,25.0) circle (4.0pt);
\node[black,anchor=center] at (27.1375,25.0) {\footnotesize 1};\draw[orange,fill=white] (26.0125,26.215) circle (4.0pt);
\node[black,anchor=center] at (26.0125,26.215) {\footnotesize 3};\draw[black] (27.925,25.315) -- (27.3625,26.215);
\draw[cyan,fill=white] (27.64375,25.765) circle (4.0pt);
\node[black,anchor=center] at (27.64375,25.765) {\footnotesize 0};\draw[black] (28.4875,25.0) -- (27.925,25.315);
\draw[cyan,fill=white] (28.20625,25.1575) circle (4.0pt);
\node[black,anchor=center] at (28.20625,25.1575) {\footnotesize 2};\draw[black] (27.3625,25.09) -- (27.925,25.315);
\draw[cyan,fill=white] (27.64375,25.2025) circle (4.0pt);
\node[black,anchor=center] at (27.64375,25.2025) {\footnotesize 4};\draw[black] (28.4875,25.0) -- (27.3625,26.215);
\draw[cyan,fill=white] (27.925,25.6075) circle (4.0pt);
\node[black,anchor=center] at (27.925,25.6075) {\footnotesize 1};\draw[black] (27.3625,25.09) -- (27.3625,26.215);
\draw[cyan,fill=white] (27.3625,25.6525) circle (4.0pt);
\node[black,anchor=center] at (27.3625,25.6525) {\footnotesize 3};\draw[black] (27.3625,25.09) -- (28.4875,25.0);
\draw[cyan,fill=white] (27.925,25.045) circle (4.0pt);
\node[black,anchor=center] at (27.925,25.045) {\footnotesize 5};\draw[green,fill=white] (29.275,25.51) circle (4.0pt);
\node[black,anchor=center] at (29.275,25.51) {\footnotesize 0};\draw[green,fill=white] (28.9,25.54) circle (4.0pt);
\node[black,anchor=center] at (28.9,25.54) {\footnotesize 1};\draw[green,fill=white] (29.275,25.135) circle (4.0pt);
\node[black,anchor=center] at (29.275,25.135) {\footnotesize 3};\draw[black] (29.275,25.315) -- (28.7125,26.215);
\draw[black] (29.8375,25.0) -- (29.275,25.315);
\draw[black] (28.7125,25.09) -- (29.275,25.315);
\draw[green,fill=white] (29.0875,25.435) circle (4.0pt);
\node[black,anchor=center] at (29.0875,25.435) {\footnotesize 2};\draw[black] (29.8375,25.0) -- (28.7125,26.215);
\draw[black] (28.7125,25.09) -- (28.7125,26.215);
\draw[black] (28.7125,25.09) -- (29.8375,25.0);
\draw[black] (30.625,25.315) -- (30.0625,26.215);
\draw[black] (31.1875,25.0) -- (30.625,25.315);
\draw[black] (30.0625,25.09) -- (30.625,25.315);
\draw[magenta,fill=white] (30.484375,25.405) circle (4.0pt);
\node[black,anchor=center] at (30.484375,25.405) {\footnotesize 0};\draw[black] (31.1875,25.0) -- (30.0625,26.215);
\draw[black] (30.0625,25.09) -- (30.0625,26.215);
\draw[black] (30.0625,25.09) -- (31.1875,25.0);
\end{tikzpicture}
\caption{The numbering of the sub-entities of a tetrahedral reference cell.}
\label{fig:tet-reference}
\end{figure}

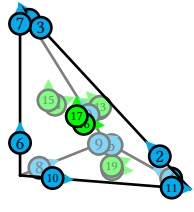
\begin{figure}
\begin{tikzpicture}[x=1cm,y=1cm,line cap=round,line join=round,line width=1pt]
\fill[white,opacity=0.5] (27.279999602501192,25.002904458356085) -- (26.154999602501192,25.632904458356084) -- (25.029999602501192,27.432904458356084) -- cycle;\fill[white,opacity=0.5] (25.029999602501192,25.182904458356084) -- (26.154999602501192,25.632904458356084) -- (25.029999602501192,27.432904458356084) -- cycle;\fill[white,opacity=0.5] (25.029999602501192,25.182904458356084) -- (27.279999602501192,25.002904458356085) -- (26.154999602501192,25.632904458356084) -- cycle;\draw[black] (26.154999602501192,25.632904458356084) -- (25.029999602501192,27.432904458356084);
\draw[black] (27.279999602501192,25.002904458356085) -- (26.154999602501192,25.632904458356084);
\draw[black] (25.029999602501192,25.182904458356084) -- (26.154999602501192,25.632904458356084);
\draw[-stealth,green] (26.073202163505908,26.068711024193444) -- (25.827809846520047,26.206130721705524);
\draw[black,fill=green] (26.073202163505908,26.068711024193444) circle (4.0pt);
\node[black,anchor=center] at (26.073202163505908,26.068711024193444) {\tiny 12};\draw[-stealth,green] (26.091305099548226,26.091694521545286) -- (25.900221590689327,26.298064711112897);
\draw[black,fill=green] (26.091305099548226,26.091694521545286) circle (4.0pt);
\node[black,anchor=center] at (26.091305099548226,26.091694521545286) {\tiny 13};\draw[-stealth,green] (25.492044292271686,26.117722334264283) -- (25.753178361583164,26.222175961988874);
\draw[black,fill=green] (25.492044292271686,26.117722334264283) circle (4.0pt);
\node[black,anchor=center] at (25.492044292271686,26.117722334264283) {\tiny 14};\draw[-stealth,green] (25.404999602501192,26.176654458356083) -- (25.404999602501192,26.457904458356083);
\draw[black,fill=green] (25.404999602501192,26.176654458356083) circle (4.0pt);
\node[black,anchor=center] at (25.404999602501192,26.176654458356083) {\tiny 15};\draw[-stealth,green] (26.248451034863955,25.265428343767063) -- (26.52880533195224,25.243);
\draw[black,fill=green] (26.248451034863955,25.265428343767063) circle (4.0pt);
\node[black,anchor=center] at (26.248451034863955,25.265428343767063) {\tiny 18};\draw[-stealth,green] (26.242044292271686,25.30772233426428) -- (26.503178361583164,25.412175961988872);
\draw[black,fill=green] (26.242044292271686,25.30772233426428) circle (4.0pt);
\node[black,anchor=center] at (26.242044292271686,25.30772233426428) {\tiny 19};\draw[-stealth,cyan] (25.936562201875894,25.98240429935656) -- (25.7875,26.220903822357993);
\draw[black,fill=cyan] (25.936562201875894,25.98240429935656) circle (4.0pt);
\node[black,anchor=center] at (25.936562201875894,25.98240429935656) {\footnotesize 0};\draw[-stealth,cyan] (25.149062201875893,27.24240429935656) -- (25.0,27.48090382235799);
\draw[black,fill=cyan] (25.149062201875893,27.24240429935656) circle (4.0pt);
\node[black,anchor=center] at (25.149062201875893,27.24240429935656) {\footnotesize 1};\draw[-stealth,cyan] (27.029452163505905,25.143211024193445) -- (26.784059846520048,25.280630721705524);
\draw[black,fill=cyan] (27.029452163505905,25.143211024193445) circle (4.0pt);
\node[black,anchor=center] at (27.029452163505905,25.143211024193445) {\footnotesize 4};\draw[-stealth,cyan] (26.241952163505907,25.584211024193444) -- (25.996559846520046,25.721630721705523);
\draw[black,fill=cyan] (26.241952163505907,25.584211024193444) circle (4.0pt);
\node[black,anchor=center] at (26.241952163505907,25.584211024193444) {\footnotesize 5};\draw[-stealth,cyan] (25.285794292271685,25.28522233426428) -- (25.546928361583166,25.38967596198887);
\draw[black,fill=cyan] (25.285794292271685,25.28522233426428) circle (4.0pt);
\node[black,anchor=center] at (25.285794292271685,25.28522233426428) {\footnotesize 8};\draw[-stealth,cyan] (26.073294292271687,25.60022233426428) -- (26.334428361583164,25.704675961988873);
\draw[black,fill=cyan] (26.073294292271687,25.60022233426428) circle (4.0pt);
\node[black,anchor=center] at (26.073294292271687,25.60022233426428) {\footnotesize 9};\fill[white,opacity=0.5] (25.029999602501192,25.182904458356084) -- (27.279999602501192,25.002904458356085) -- (25.029999602501192,27.432904458356084) -- cycle;\draw[black] (27.279999602501192,25.002904458356085) -- (25.029999602501192,27.432904458356084);
\draw[black] (25.029999602501192,25.182904458356084) -- (25.029999602501192,27.432904458356084);
\draw[black] (25.029999602501192,25.182904458356084) -- (27.279999602501192,25.002904458356085);
\draw[-stealth,green] (25.873451034863955,25.865428343767064) -- (26.15380533195224,25.843);
\draw[black,fill=green] (25.873451034863955,25.865428343767064) circle (4.0pt);
\node[black,anchor=center] at (25.873451034863955,25.865428343767064) {\tiny 16};\draw[-stealth,green] (25.779999602501192,25.966654458356082) -- (25.779999602501192,26.247904458356082);
\draw[black,fill=green] (25.779999602501192,25.966654458356082) circle (4.0pt);
\node[black,anchor=center] at (25.779999602501192,25.966654458356082) {\tiny 17};\draw[-stealth,cyan] (26.878805099548227,25.436194521545286) -- (26.68772159068933,25.642564711112897);
\draw[black,fill=cyan] (26.878805099548227,25.436194521545286) circle (4.0pt);
\node[black,anchor=center] at (26.878805099548227,25.436194521545286) {\footnotesize 2};\draw[-stealth,cyan] (25.303805099548228,27.137194521545286) -- (25.11272159068933,27.343564711112897);
\draw[black,fill=cyan] (25.303805099548228,27.137194521545286) circle (4.0pt);
\node[black,anchor=center] at (25.303805099548228,27.137194521545286) {\footnotesize 3};\draw[-stealth,cyan] (25.029999602501192,25.614154458356083) -- (25.029999602501192,25.895404458356083);
\draw[black,fill=cyan] (25.029999602501192,25.614154458356083) circle (4.0pt);
\node[black,anchor=center] at (25.029999602501192,25.614154458356083) {\footnotesize 6};\draw[-stealth,cyan] (25.029999602501192,27.189154458356082) -- (25.029999602501192,27.470404458356082);
\draw[black,fill=cyan] (25.029999602501192,27.189154458356082) circle (4.0pt);
\node[black,anchor=center] at (25.029999602501192,27.189154458356082) {\footnotesize 7};\draw[-stealth,cyan] (25.460951034863953,25.148428343767062) -- (25.74130533195224,25.126);
\draw[black,fill=cyan] (25.460951034863953,25.148428343767062) circle (4.0pt);
\node[black,anchor=center] at (25.460951034863953,25.148428343767062) {\tiny 10};\draw[-stealth,cyan] (27.035951034863956,25.022428343767064) -- (27.31630533195224,25.0);
\draw[black,fill=cyan] (27.035951034863956,25.022428343767064) circle (4.0pt);
\node[black,anchor=center] at (27.035951034863956,25.022428343767064) {\tiny 11};\end{tikzpicture}\caption{The DOFs of a degree 2 N\'ed\'elec first kind element on a tetrahedron.}
\label{fig:tetra-dofs}
\end{figure}

\begin{example}[N\'ed\'elec degree 2 on a tetrahedron]
\label{example-def-nedelec}
A degree 2 N\'ed\'elec first kind element on a tetrahedral cell
\cite{nedelec} is defined by
\begin{itemize}
  \item $\refel = \left\{(x,y,z)\in[0,1]^3\,\middle|\,x+y+z\leqslant1\right\}$
 (where we number the sub-entities as shown in
 \cref{fig:tet-reference}),

  \item $\polyspace = \mathbb{P}_1^3 \oplus \left\{p\in
  \mathbb{P}_2^3\setminus\mathbb{P}_1^3\,\middle|\,p\cdot
  \begin{bmatrix}x\\y\\z\end{bmatrix} = 0\right\}$,
where $\mathbb{P}_k = \left\{x^iy^jz^m\,\middle|\,i,j,m \in
\{0,1,\dots,k\}\text{ and } i + j + m \leqslant k \right\}$,
  \item $\dualbasis=\{\hat{l}_0,\dots,\hat{l}_{19}\}$.
The
functionals $\hat{l}_0$ and $\hat{l}_1$ are integral moments of
tangential components against two linear functions on edge 0, defined by
\begin{align*}
  \hat{l}_{0} &: v \mapsto \int_0^1v(0, 1-t, t) \cdot
    \begin{bmatrix}0\\t-1\\1-t\end{bmatrix}\,\mathrm{d}t,&
  \hat{l}_{1} &: v \mapsto \int_0^1v(0, 1-t, t) \cdot
    \begin{bmatrix}0\\-t\\t\end{bmatrix}\,\mathrm{d}t.
\end{align*}
The functionals $\hat{l}_2$ to $\hat{l}_{11}$ are the corresponding
integral moments for edges 1 to 5. The functionals $\hat{l}_{18}$ and
$\hat{l}_{19}$ are integral moments of the two tangential components
against a constant function on face 3, defined by
\begin{align*}
\hat{l}_{18}&:v\mapsto\int_0^1\int_0^{1-t}v(s, t, 0 )\cdot
    \begin{bmatrix}1\\0\\0\end{bmatrix}\,\mathrm{d}s\,\mathrm{d}t,&
\hat{l}_{19}&:v\mapsto\int_0^1\int_0^{1-t}v(s, t, 0 )\cdot
    \begin{bmatrix}0\\1\\0\end{bmatrix}\,\mathrm{d}s\,\mathrm{d}t.
\end{align*}
The functionals $\hat{l}_{12}$ to $\hat{l}_{17}$ are the corresponding
functionals on faces 0 to 2.
\end{itemize}
In Basix, these functionals are implemented with a set of quadrature
points and a tensor whose entries are the product of a quadrature
weight, the value of the linear function at the quadrature point, and a
component of the normal vector. A representation of the DOFs of this
element is shown in \cref{fig:tetra-dofs}.

The push-forward for this element is the covariant Piola map. If
$g:\refel \to \mathbb{R}^d$ maps points from the reference cell to a
physical cell, then the covariant Piola map
$\mathcal{F}_g^\textup{curl}$ is defined by
$$
  \mathcal{F}_{g}^\textup{curl}(\hat{\psi})
    = \mat{J}_{g}\invtranspose\hat{\psi} \circ g^{-1},
$$
where $\mat{J}_g$ is the Jacobian of $g$ and $\hat{\psi}$ is a
vector-valued function on the reference cell. The covariant Piola map
preserves the tangential components of functions when they are pushed
forward.
\end{example}

\section{Degree-of-freedom transformations}
\label{sec:dof-transformations}

\begin{figure}
\begin{tikzpicture}[line cap=round,line join=round,line width=1pt]
\footnotesize
\draw (0,0) -- (0,2);
\draw (1,0) -- (1,2);
\draw[dashed] (0,2/3) -- (1,4/3) (0,4/3) -- (1,2/3);
\draw[fill=cyan] (0,2/3) circle (1.5mm) node {0};
\draw[fill=cyan] (0,4/3) circle (1.5mm) node {1};
\draw[fill=cyan] (1,2/3) circle (1.5mm) node {1};
\draw[fill=cyan] (1,4/3) circle (1.5mm) node {0};
\end{tikzpicture}
\hspace{16mm}
\begin{tikzpicture}[line cap=round,line join=round,line width=1pt]
\footnotesize
\draw[fill=white] (0,0) -- (1.5,0.5) -- (0,2) -- cycle;
\draw[->,orange] (0.7*0.5+0.15,0.7*2.5/3+0.05) -- +(0.5,0.5/3);
\draw[->,orange] (0.7*0.5,0.7*2.5/3+0.2) -- +(0,2/3);
\draw[fill=orange] (0.7*0.5+0.15,0.7*2.5/3+0.05) circle (1.5mm) node {0};
\draw[fill=orange] (0.7*0.5,0.7*2.5/3+0.2) circle (1.5mm) node {1};

\begin{scope}[shift={(1.3,0)}]
\draw[fill=white] (0,0) -- (1.5,0.5) -- (0,2) -- cycle;
\draw[->,orange] (0.3*1.5+0.7*0.5-0.15,0.3*0.5+0.7*2.5/3+0.15) -- +(-0.5,0.5);
\draw[->,orange] (0.3*1.5+0.7*0.5-0.15,0.3*0.5+0.7*2.5/3-0.05) -- +(-0.5,-0.5/3);
\draw[fill=orange] (0.3*1.5+0.7*0.5-0.15,0.3*0.5+0.7*2.5/3+0.15) circle (1.5mm) node {0};
\draw[fill=orange] (0.3*1.5+0.7*0.5-0.15,0.3*0.5+0.7*2.5/3-0.05) circle (1.5mm) node {1};
\end{scope}
\end{tikzpicture}
\caption{The two vertical lines on the left represent the same edge from
the point of view of two neighbouring cells, and include the DOFs of a
degree 3 Lagrange element that are associated with the edge. If the
cells do not agree on the orientation of the edge, then the wrong DOFs
will be assigned the same global DOF number. This would lead to an
incorrect basis function, as shown in \cref{fig:mismatched_edge}. The
two triangles on the right represent the same face from the point of
view of two neighbouring cells, and include the DOFs of a degree 2
N\'ed\'elec first kind element that are associated with the face. If the
cells do not agree on the orientation of the face, then the directions
of the normals to the face will not agree. For higher degree elements,
there will also be an incorrect combination of local DOFs.}
\label{fig:mismatch}
\end{figure}
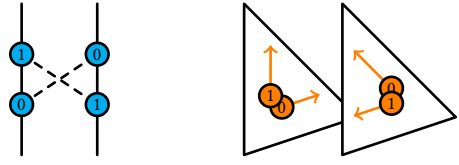

In this section, we examine how degree-of-freedom transformations can be
used to account for differences between the orientations of sub-entities
on physical cells compared to the reference cell. Not accounting for
such differences would lead to mismatches in the basis functions
associated with shared sub-entities, as illustrated in
\cref{fig:mismatch}.

\subsection{Degree-of-freedom transformation matrix}
\label{sec:dof-transformations-matrix}

As proposed in~\cite{2022-dofs}, we will define a
DOF transformation matrix $\mat{T}$ for each cell such that the basis
functions $\phi_0, \dots, \phi_{n-1}$ on a physical cell whose entities
have been reflected and rotated to match the orientation of the
sub-entities of its neighbouring cells are given by
\begin{equation}
  \label{define-T}
  \hat{\vec{\phi}}_g \coloneq
  \begin{bmatrix}
    \mathcal{F}_g(\hat{\phi}_0)\\ \vdots \\ \mathcal{F}_g(\hat{\phi}_{n-1})
  \end{bmatrix} = \mat{T}
  \begin{bmatrix}
    \phi_0\\\vdots\\\phi_{n-1}
  \end{bmatrix} \eqcolon
 \mat{T} \vec{\phi},
\end{equation}
where $\mathcal{F}_g$ is the push forward corresponding to the geometry
map $g$ from the reference cell to the physical cell. The role of the
matrix $\mat{T}$ is to account for differences between the orientation
of the sub-entities of the reference cell and the corresponding
sub-entities of the physical cell.

A finite element function $f$ at some point within a cell can be
evaluated via
\begin{equation}
  f
  = \vec{\phi}^{T} \vec{c}
  = \hat{\vec{\phi}}_g^{T} \hat{\vec{c}},
  \label{eq:transform_dot}
\end{equation}
where $\vec{\phi}$ and $\vec{c}$ are vectors containing the basis
functions and degrees-of-freedom, respectively, relative to the physical
cell ordering, and $\hat{\vec{\phi}}_g$ and $\hat{\vec{c}}$ are the basis
functions and degrees-of-freedom, respectively, relative to the
reference cell ordering. Using \cref{define-T}, we see that
$\hat{\vec{\phi}}_g^{T} \hat{\vec{c}} = (\mat{T} \vec{\phi})^{T}
\hat{\vec{c}} = \vec{\phi}^{T} (\mat{T}^{T} \hat{\vec{c}})$, which shows
that $\vec{c} = \mat{T}^{T} \hat{\vec{c}}$. It follows trivially that
$\vec{\phi} = \mat{T}^{-1} \hat{\vec{\phi}}_g$ and $\hat{\vec{c}} =
\mat{T}^{-T} \vec{c}$. If a finite element matrix $\hat{\mat{A}} \in
\mathbb{C}^{n_{0} \times n_{1}}$ for a cell is computed following the
reference cell ordering of basis functions, it is straightforward to
show the element matrix for the physical cell ordering is given by
$\mat{A} = \mat{T}_{0}^{T} \hat{\mat{A}} \mat{T}_{1}^{-T}$, where
$\mat{T}_{0}$ and $\mat{T}_{1}$ are the transformation matrices for the
trial and test function spaces, respectively. The objective is to
compute $\mat{T}$, with which we can compute element matrices and
vectors using the reference element ordering, and transform these to a
consistent global ordering.

\subsection{Cell sub-entity rotation and reflection and geometric mapping}
\label{sec:sub-entity-transforming}

We begin by defining transformations of cell sub-entities that allow an
entity to be transformed to any orientation. Orientation of a vertex
(dimension~0) is trivial, and DOFs associated with the interior of a
cell are not shared by more than one cell, hence orientation of a cell
does not require consideration. We focus here on edges (dimension~1) and
faces (dimension~2).

The orientation of an interval (edge) can only be changed by a
reflection. We can therefore express any re-orientation of the interval
as $\left(G^{\interval}_{\reflection}\right)^\alpha$, where
$G^{\interval}_{\reflection} : \mathbb{R} \rightarrow \mathbb{R}$ is a
reflection of the interval and $\alpha \in \{0, 1\}$.
Two-dimensional sub-entities are polygons. The orientation of a polygon
can be changed by rotation and/or reflection: if we define a rotation
$G^{\twodimentity}_{\rotation} : \mathbb{R}^{2} \rightarrow
\mathbb{R}^{2}$ and a reflection $G^{\twodimentity}_{\reflection} :
\mathbb{R}^{2} \rightarrow \mathbb{R}^{2}$, then we can write any
re-orientation of the polygon as
$\left(G^{\twodimentity}_{\reflection}\right)^\alpha \circ
\left(G^{\twodimentity}_{\rotation}\right)^\beta$, where $\alpha \in
\{0, 1\}$ and $\beta \in \{0, 1, \dots, \text{number of vertices} -
1\}$. For $G^{\twodimentity}_{\rotation}$, we pick a clockwise rotation
by one vertex, and for $G^{\twodimentity}_{\reflection}$ we pick a
reflection in the line $y = x$. For standard finite element cells
(tetrahedra, hexahedra, {\prismname}s, and {\pyramidname}s),
two-dimensional sub-entities are triangles or quadrilaterals. The values
of $G^{\twodimentity}_{\rotation}$ and $G^{\twodimentity}_{\reflection}$
that we use for these are given in \cref{table:cells}.
In the language of group theory, the transformations
$G^{\interval}_{\reflection}$ (for an edge) and
$G^{\twodimentity}_{\rotation}$ and $G^{\twodimentity}_{\reflection}$
(for a face) are generators of the symmetry group of the sub-entity.

\begin{table}
  \def\arraystretch{1.8}
  \begin{tabular}{m{2.5cm}|m{3cm}|m{4cm}}
  \textbf{sub-entity type}&\centering\textbf{reference sub-entity}&\textbf{generators}\\\hline
  interval&\centering\begin{tikzpicture}[line width=1pt,line cap=round,line join=round]
  \footnotesize
  \draw (0,0) -- (1,0);
  \node[anchor=south] at (0,0) {0};
  \node[anchor=south] at (1,0) {1};
  \end{tikzpicture}&$G^{\interval}_{\reflection}:x\mapsto1-x$\\
  triangle&\centering\begin{tikzpicture}[line width=1pt,line cap=round,line join=round]
  \footnotesize
  \draw (0,0) -- (1,0) -- (0,1) -- cycle;
  \node[anchor=north] at (0,0) {(0,0)};
  \node[anchor=north] at (1,0) {(1,0)};
  \node[anchor=south] at (0,1) {(0,1)};
  \end{tikzpicture}&$G^{\rtriangle}_{\rotation}:(x,y)\mapsto(y, 1-x-y)$\newline$G^{\rtriangle}_{\reflection}:(x,y)\mapsto(y,x)$\\
  quadrilateral&\centering\begin{tikzpicture}[line width=1pt,line cap=round,line join=round]
  \footnotesize
  \draw (0,0) -- (1,0) -- (1, 1) -- (0,1) -- cycle;
  \node[anchor=north] at (0,0) {(0,0)};
  \node[anchor=north] at (1,0) {(1,0)};
  \node[anchor=south] at (0,1) {(0,1)};
  \node[anchor=south] at (1,1) {(1,1)};
  \end{tikzpicture}&$G^{\quadrilateral}_{\rotation}:(x,y)\mapsto(y, 1-x)$\newline$G^{\quadrilateral}_{\reflection}:(x,y)\mapsto(y,x)$\\
  \end{tabular}
  \vspace{2mm}
  \caption{The cell sub-entities and the transformations we use to
  generate their symmetry groups. We represent each sub-entity in this
  table using the reference used in Symfem and Basix, although this is
  arbitrary. If a different reference (for example, some use the
  reference quadrilateral with vertices at $(\pm1,\pm1)$) is used, the
  generator functions must be updated.}
  \label{table:cells}
\end{table}

\newcommand{\mapsymbol}{H}

We also define $\hat{E} \subset \mathbb{R}^{d_E}$ to be the reference
cell with the same cell type as the sub-entity $E \subset
\mathbb{R}^{d_R}$; where $0 < d_E < d_R$, and introduce an affine map
$g_E : \mathbb{R}^{d_E} \to \mathbb{R}^{d_{\refel}}$ that maps each
vertex of $\hat{E}$ to a corresponding vertex of $E$. If $E$ is an edge, let
$\mapsymbol^E_{\reflection}$ be any linear invertible map such that
\begin{subequations}\label{define_Hmaps}
  \begin{align}
    \left.\mapsymbol^E_{\reflection}\right|_E
      &= g_E \circ G^{\interval}_{\reflection} \circ g^{-1}_E,
  \label{define_Hmaps_edge}
  \end{align}
  where $\left. \mapsymbol^E_{\reflection} \right|_E$ is the restriction
  of $\mapsymbol^E_{\reflection}$ to $E$. If $E$ is a face, let
  $\mapsymbol^E_{\reflection}$ and $\mapsymbol^E_{\rotation}$
  be any linear invertible maps such that
  \begin{align}
    \left.\mapsymbol^E_{\reflection}\right|_E
    &\coloneq g_E \circ G^{\twodimentity}_{\reflection} \circ g^{-1}_E,
    &\left.\mapsymbol^E_{\rotation}\right|_E
      &\coloneq g_E \circ G^{\twodimentity}_{\rotation} \circ g^{-1}_E.
  \end{align}
\end{subequations}

\begin{figure}
\begin{tikzpicture}[x=2cm,y=2cm,line cap=round,line join=round,line width=1pt]
\footnotesize
\newcommand{\zeromarker}[2][0]{\begin{scope}[shift={#2},rotate=#1]
\draw[orange,fill=white]
    ({sin(0)*2.3pt},{cos(0)*2.3pt}) --
    ({sin(120)*2.3pt},{cos(120)*2.3pt}) --
    ({sin(240)*2.3pt},{cos(240)*2.3pt}) --
    cycle;
\end{scope}}
\newcommand{\onemarker}[1]{\begin{scope}[shift={#1}]
\draw[orange,fill=white] (-1.7pt,-1.7pt) rectangle (1.7pt,1.7pt);
\end{scope}}
\draw[->,line width=0.5pt] (0,0) -- (1.2,0) node[anchor=west] {$x$};
\draw[->,line width=0.5pt] (0,0) -- (0,1.2) node[anchor=south] {$y$};
\draw (0,0) -- (1,0) -- (0,1);
\draw[orange] (0,0) -- (0,1);
\node[anchor=east] at (0,0) {$g_E(0)=(0,0)$};
\node[anchor=north] at (1,0) {$(1,0)$};
\node[anchor=east] at (0,1) {$g_E(1)=(0,1)$};
\node at (0.33333,0.33333) {$\refel$};
\node[orange,anchor=east] at (0,0.5) {$E$};
\zeromarker[180]{(0,0)}
\onemarker{(0,1)}

\begin{scope}[shift={(2.5,0)}]
\draw[->,line width=0.5pt] (0,0) -- (1.2,0) node[anchor=west] {$x$};
\draw[->,line width=0.5pt] (0,0) -- (0,1.2) node[anchor=south] {$y$};
\draw (0,0) -- (1,1) -- (0,1);
\draw[orange] (0,0) -- (0,1);
\node[anchor=north] at (0,0) {$\mapsymbol^E_{\reflection}((0,1))$};
\node[anchor=east] at (0,1) {$\mapsymbol^E_{\reflection}((0,0))$};
\node[anchor=west] at (1,1) {$\mapsymbol^E_{\reflection}((1,0))$};
\node[orange,anchor=east] at (0,0.4) {$\mapsymbol^E_{\reflection}(E)$};
\node at (0.33333,0.66666) {$\mapsymbol^E_{\reflection}(R)$};
\zeromarker{(0,1)}
\onemarker{(0,0)}
\end{scope}

\begin{scope}[shift={(-2.5,0)}]
\draw[->,line width=0.5pt] (0,0) -- (1.2,0) node[anchor=west] {$x$};
\draw[orange] (0,0) -- (1,0);
\node[anchor=north] at (0.5,0) {$\hat{E}$};
\node[anchor=north] at (0,0) {$0$};
\node[anchor=north] at (1,0) {$1$};
\zeromarker[90]{(0,0)}
\onemarker{(1,0)}
\end{scope}

\draw[cyan,->,shorten >=10pt,shorten <= 5pt] (-2,0) -- (0,0.5);
\node[cyan,fill=white] at (-1,0.25) {$g_E$};

\draw[cyan,->,shorten >=20pt,shorten <=20pt] (0.33333,0.33333) -- (2.83333,0.66666);
\node[cyan,fill=white] at (1.58333,0.5) {$\mapsymbol^E_{\reflection}$};

\end{tikzpicture}
\caption[An example showing the notation introduced.]{An example showing
the notation introduced. In this example, $\refel$ is the reference
triangle, $E$ is an edge, $d_\refel = 2$, $d_E = 1$,
$g_E:x\mapsto(0,x)$,
$G^{\interval}_{\reflection} : x \mapsto 1 - x$, and
$\mapsymbol^E_{\reflection} : (x,y) \mapsto (x,1-y)$. As stated in
\cref{define_Hmaps_edge}, $\left.\mapsymbol^E_{\reflection}\right|_E=g_E
\circ G^{\interval}_{\reflection} \circ g^{-1}_E$.
}
\label{fig:computing_example}
\end{figure}
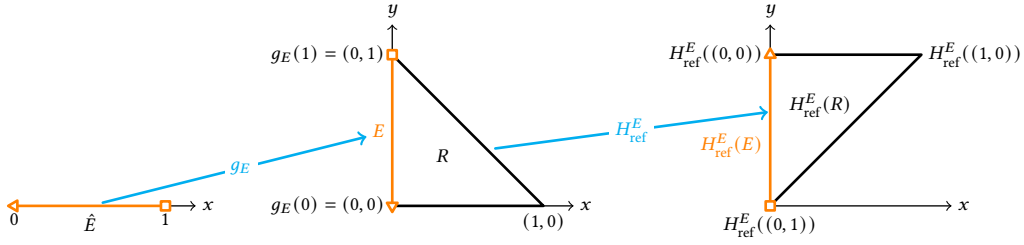

When applied to a point on the sub-entity $E$, the operators $H^E_{\reflection}$ and $H^E_{\rotation}$
(i) map the point to the reference sub-entity $\hat{E}$ (by
application of $g^{-1}_E$), (ii) rotate or reflect the sub-entity
reference cell ($G^{\twodimentity}_{\reflection}$ or $G^{\twodimentity}_{\rotation}$), and then
(iii) map the point back to the sub-entity $E$ ($g_E$).
\Cref{fig:computing_example} illustrates this for a case where where
$\refel$ is a triangle.

\subsection{Basis function and degree-of-freedom transformations}
\label{sec:split_by_entity}

The DOF transformation matrix $\mat{T}$ in \cref{define-T} represents
the effect of reorienting the sub-entities of a cell on the basis
functions. Due to \cref{assume-me2}, changing the orientation of a sub-entity of a cell will only
affect the DOFs associated with that sub-entity, hence we can consider
the DOFs associated with each sub-entity separately. The
entity-by-entity numbering that we use (see \cref{entity-by-entity})
means that $\mat{T}$ for a cell will be block-diagonal with a block for
each cell sub-entity. For example, for an element on a triangular cell,
$$
\mat{T} = \begin{bmatrix}
\mat{T}_{00}&0&0&0&0&0&0\\
0&\mat{T}_{01}&0&0&0&0&0\\
0&0&\mat{T}_{02}&0&0&0&0\\
0&0&0&\mat{T}_{10}&0&0&0\\
0&0&0&0&\mat{T}_{11}&0&0\\
0&0&0&0&0&\mat{T}_{12}&0\\
0&0&0&0&0&0&\mat{T}_{20}\\
\end{bmatrix},
$$
where $\mat{T}_{ij} \in \mathbb{R}^{n_{ij} \times n_{ij}}$ is the
block for the $j$th entity of dimension $i$, and $n_{ij}$ is the number
of DOFs associated with the $j$th entity of dimension~$i$.

As noted in \cref{sec:sub-entity-transforming}, we do not need to apply
DOF transformations to basis functions associated with vertices or the
interior of the cell, so the blocks of $\mat{T}$ for these entities are
the identity. In this triangle example, this means that $\mat{T}_{00}$,
$\mat{T}_{01}$, $\mat{T}_{02}$ and $\mat{T}_{20}$ are all identity
matrices, and
\begin{equation}\label{eq:blocks}
\mat{T} = \begin{bmatrix}
\mat{I} & 0 & 0 & 0 & 0\\
0 & \mat{T}_{10} & 0 & 0 &0\\
0 & 0 & \mat{T}_{11} & 0 &0\\
0 & 0 &0 &\mat{T}_{12} & 0\\
0 & 0 &0 & 0 &\mat{I}\\
\end{bmatrix}.
\end{equation}
The task now is to determine an expression for the non-trivial blocks
$\mat{T}_{ij}$.

\subsection{Computing the base transformations}

Each block $\mat{T}_{ij}$ of the DOF transformation matrix $\mat{T}$
describes the effect of applying a transformation to a sub-entity of the
cell. Let $E$ be the sub-entity of the cell corresponding to the block
$\mat{T}_{ij}$, and let $\hat{E}$ be the reference cell with the same
cell type as $E$. For simplicity, we assume in this section that $E$ is
a 2-dimensional sub-entity. What follows can be adapted trivially to
1-dimensional sub-entities by discarding the rotation
$G^{\twodimentity}_{\rotation}$ and replacing
$G^{\twodimentity}_{\reflection}$ with $G^{\interval}_{\reflection}$.

As described in \cref{sec:sub-entity-transforming}, any transformation
of the reference $\hat{E}$ can be written as a combination of the
generators $G^{\twodimentity}_{\rotation}$ and
$G^{\twodimentity}_{\reflection}$. In this section, we will define
\emph{base transformation} matrices for each sub-entity $E$ that
describe the effect on the basis functions of applying the generators
to~$E$.

Let $\hat{l}_a, \dots, \hat{l}_b$ and $\hat{\phi}_a, \dots,
\hat{\phi}_b$ be the functionals and basis functions associated with
$E$. Let $\mapsymbol^E_{\rotation}$ and $\mapsymbol^E_{\reflection}$ be
defined as in~\cref{define_Hmaps} with the additional assumptions that
$\mapsymbol^{E}_{\rotation}(\refel) = \refel$ and
$\mapsymbol^{E}_{\reflection}(\refel) = \refel$. We treat these maps as
geometry maps, and introduce the corresponding push-forward maps
$\mathcal{F}_{\mapsymbol^E_{\rotation}}$ and
$\mathcal{F}_{\mapsymbol^E_{\reflection}}$. We proceed for the
$\mathcal{F}_{\mapsymbol^E_{\reflection}}$ case; the same steps can be
followed for $\mathcal{F}_{\mapsymbol^E_{\rotation}}$. Using
\cref{assume-me3}, we know that
\begin{equation}\label{spans-same}
\polyspace_E
\coloneq \operatorname{span}\left(\left\{\hat{\phi}_a,
  \dots, \hat{\phi}_b \right\}\right)
=
\operatorname{span}\left(\left\{\mathcal{F}_{\mapsymbol^E_{\reflection}}(\hat{\phi}_a),
  \dots,
  \mathcal{F}_{\mapsymbol^E_{\reflection}}(\hat{\phi}_b)\right\}\right).
\end{equation}
Using \cref{assume-me2,ciarlet-def}, we have, for $a \leqslant i, j
\leqslant b$,
\begin{equation}\label{deltaij}
  \hat{l}_j\left(\hat{\phi}_i\right) = \delta_{ij}.
\end{equation}
From \cref{deltaij} it follows that, for any $f\in\polyspace_E$, we can
express $f$ as
$$
f = \sum_{j=a}^b\hat{l}_j(f)\hat{\phi}_j.
$$
From \cref{spans-same}, we see that for $a \leqslant i \leqslant b$, the
function $f = \mathcal{F}_{\mapsymbol^E_{\reflection}}(\hat{\phi}_i)$ is
in $\polyspace_E$, and so
$$
\mathcal{F}_{\mapsymbol^E_{\reflection}}(\hat{\phi}_i)
  = \sum_{j=a}^b\hat{l}_j\left( \mathcal{F}_{\mapsymbol^E_{\reflection}}(\hat{\phi}_i) \right)\hat{\phi}_j.
$$
In matrix form, this is
\begin{equation}
\begin{bmatrix}
  \mathcal{F}_{\mapsymbol^E_{\reflection}}(\hat{\phi}_a)\\
  \vdots\\
  \mathcal{F}_{\mapsymbol^E_{\reflection}}(\hat{\phi}_b)
  \end{bmatrix}
  =
  \underbrace{
  \begin{bmatrix}
  \hat{l}_a\left(\mathcal{F}_{\mapsymbol^E_{\reflection}}(\hat{\phi}_a)\right)
    & \dots
    & \hat{l}_b\left( \mathcal{F}_{\mapsymbol^E_{\reflection}}(\hat{\phi}_a) \right)\\
  \vdots&\ddots&\vdots
  \\
  \hat{l}_a\left(\mathcal{F}_{\mapsymbol^E_{\reflection}}(\hat{\phi}_b)\right)
    & \dots
    &\hat{l}_b\left( \mathcal{F}_{\mapsymbol^E_{\reflection}}(\hat{\phi}_b) \right)
  \end{bmatrix}
  }_{\mat{B}^E_{\reflection}}
  \begin{bmatrix}
  \hat{\phi}_a\\
  \vdots\\
  \hat{\phi}_b
\end{bmatrix},
\label{eq:B_ref}
\end{equation}
where we refer to $\mat{B}^E_{\reflection}$ as the `reflection base
transformation matrix'. Similarly for the `rotation base transformation
matrix', $\mat{B}^E_{\rotation}$:
\begin{equation}
\mat{B}^E_{\rotation} \coloneq\begin{bmatrix}
\hat{l}_a\left( \mathcal{F}_{\mapsymbol^E_{\rotation}}(\hat{\phi}_a) \right)
  &\dots &\hat{l}_b\left( \mathcal{F}_{\mapsymbol^E_{\rotation}}(\hat{\phi}_a) \right)
\\
\vdots & \ddots & \vdots
\\
\hat{l}_a\left( \mathcal{F}_{\mapsymbol^E_{\rotation}}(\hat{\phi}_b) \right)
&\dots &\hat{l}_b\left( \mathcal{F}_{\mapsymbol^E_{\rotation}}(\hat{\phi}_b) \right)
\end{bmatrix}.
\label{eq:B_rot}
\end{equation}

Note that due to \cref{assume-me2}, the functionals
$\hat{l}_a,\dots,\hat{l}_b$ only depend on the values of
$\mathcal{F}_{\mapsymbol^E_{\rotation}}(\hat{\phi}_b)$ on the sub-entity
$E$. Therefore as long as the properties in \cref{define_Hmaps} hold, we
still obtain the same matrices in \cref{eq:B_ref,eq:B_rot} if we drop
the assumptions that $\mapsymbol^{E}_{\rotation}(\refel) = \refel$ and
$\mapsymbol^{E}_{\reflection}(\refel) = \refel$, and so we will not
require these assumptions later when we define algorithms for computing
$\mat{B}^E_{\reflection}$ and $\mat{B}^E_{\rotation}$.

We now show that the two matrices $\mat{B}^E_{\reflection}$ and
$\mat{B}^E_{\rotation}$ can be used to correct for differences in
orientation by considering the effect of a push forward on the basis
functions associated with $E$. Let $g : \refel \to
\mathbb{R}^{d_\refel}$ be defined by
\begin{equation}\label{def:g=B^bB^b}
g =
  \left(\mapsymbol^{E}_{\rotation}\right)^\beta
  \circ
  \left(\mapsymbol^{E}_{\reflection}\right)^\alpha,
\end{equation}
for some $\alpha, \beta\in\mathbb{N}$. \Cref{def:g=B^bB^b} implies that
$g$ is a bijective linear map such that $g(\refel)=\refel$, $g(E) = E$,
and applying $g$ to $E$ has the same effect as applying
$
  \left(G^{\twodimentity}_{\rotation}\right)^\beta
  \circ
  \left(G^{\twodimentity}_{\reflection}\right)^\alpha
$ to $\hat{E}$. Let
$\mathcal{F}_{g}$ be the push-forward map corresponding to $g$. We
define a set of modified functionals
$\modifiedfunctional_0,\dots,\modifiedfunctional_{n-1}$ by
\begin{equation}
  \modifiedfunctional_i
  \coloneq
  \begin{cases}
    \hat{l}_i\circ\left(\mathcal{F}_g\right)^{-1}
      &a\leqslant i\leqslant b\text{ (i.e.~if $\hat{l}_i$ is associated with $E$)},\\
    \hat{l}_i &\text{otherwise}.
  \end{cases}
  \label{eqn:modified_functionals}
\end{equation}
These are the functionals that would be used to define our element if
$E$ had been differently oriented on the reference cell, as applying a
modified functional on the transformed sub-entity is equivalent to
undoing the transformation with $\left(\mathcal{F}_g\right)^{-1}$ then
applying one of our original functionals.

\Cref{def:g=B^bB^b} implies that, for $a \leqslant i \leqslant b$,
\begin{equation}\label{step1}
  \mathcal{F}_{g}\left(
  \hat{\phi}_i\right)
  =
  \left[
    \left(\mathcal{F}_{\mapsymbol^E_{\rotation}}\right)^\beta
    \circ
    \left(\mathcal{F}_{\mapsymbol^E_{\reflection}}\right)^\alpha
  \right] \left(\hat{\phi}_i \right).
\end{equation}
Using \cref{eq:B_ref} and the linearity of
$\mathcal{F}_{\mapsymbol^E_{\reflection}}$, we see that
\begin{align*}
  \left[
  \left(\mathcal{F}_{\mapsymbol^E_{\rotation}}\right)^\beta
  \circ
  \left(\mathcal{F}_{\mapsymbol^E_{\reflection}}\right)^\alpha
  \right]\left(
  \begin{bmatrix}
  \hat{\phi}_a\\\vdots\\\hat{\phi}_b
  \end{bmatrix}
  \right)
  &=
  \left[
  \left(\mathcal{F}_{\mapsymbol^E_{\rotation}}\right)^\beta
  \circ
  \left(\mathcal{F}_{\mapsymbol^E_{\reflection}}\right)^{\alpha-1}
  \right]\left(
  \mat{B}^E_{\reflection}
  \begin{bmatrix}
  \hat{\phi}_a\\\vdots\\\hat{\phi}_b
  \end{bmatrix}
  \right)\\
  &=
  \mat{B}^E_{\reflection}
  \left[
  \left(\mathcal{F}_{\mapsymbol^E_{\rotation}}\right)^\beta
  \circ
  \left(\mathcal{F}_{\mapsymbol^E_{\reflection}}\right)^{\alpha-1}
  \right]\left(
  \begin{bmatrix}
  \hat{\phi}_a\\\vdots\\\hat{\phi}_b
  \end{bmatrix}
  \right),
\end{align*}
where $\mathcal{F}(\vec{a})$ denotes the application of $\mathcal{F}$ to
each entry in a vector of functions $\vec{a}$. Applying the same step
for each $\mathcal{F}_{\mapsymbol^E_{\rotation}}$ and
$\mathcal{F}_{\mapsymbol^E_{\reflection}}$ leads to
\begin{equation}
  \left[
  \left(\mathcal{F}_{\mapsymbol^E_{\rotation}}\right)^\beta
  \circ
  \left(\mathcal{F}_{\mapsymbol^E_{\reflection}}\right)^\alpha
  \right]\left(
  \begin{bmatrix}
  \hat{\phi}_a\\\vdots\\\hat{\phi}_b
  \end{bmatrix}
  \right)=
  \left(\mat{B}^E_{\reflection}\right)^\alpha
  \left(\mat{B}^E_{\rotation}\right)^\beta
  \begin{bmatrix}
    \hat{\phi}_a
    \\\vdots\\
    \hat{\phi}_b
  \end{bmatrix}
  \eqcolon
  \begin{bmatrix}
  \modifiedbasis_a
  \\
  \vdots
  \\
  \modifiedbasis_b
  \end{bmatrix}.
\label{step2}
\end{equation}
Combining \cref{step1,step2}, we see that for $a\leqslant i\leqslant b$,
\begin{equation}\label{the-conclusion}
  \mathcal{F}_{g}(\hat{\phi}_i)
  =
  \modifiedbasis_i.
\end{equation}
From \cref{the-conclusion,spans-same}, we see that the functions
$\{\modifiedbasis_a, \dots, \modifiedbasis_b\}$ span the same space as
$\{\hat{\phi}_a, \dots, \hat{\phi}_b\}$, and so the functions
$\{\modifiedbasis_0, \dots, \modifiedbasis_{n-1}\}$ are a basis of
$\polyspace$.

Using the definition of $\modifiedfunctional_i$ (see
\cref{eqn:modified_functionals}) and \cref{the-conclusion}, we see that
for $a \leqslant i \leqslant b$ and $a \leqslant j \leqslant b$,
\begin{align}
\modifiedfunctional_i(\modifiedbasis_j)
&=
\left[\hat{l}_i\circ\left(\mathcal{F}_g\right)^{-1}\right]\left(\mathcal{F}_g(\hat{\phi}_j)\right)
\notag\\
&=
\hat{l}_i(\hat{\phi}_j)
\notag\\
&=
\delta_{ij}.\label{delta_ij}
\end{align}
From \cref{assume-me,assume-me2}, we can see that the equivalent of
\cref{delta_ij} holds for $i\not\in\{a,\dots,b\}$ and
$j\not\in\{a,\dots,b\}$, and so we conclude that the functions
$\{\modifiedbasis_0,\dots,\modifiedbasis_{n-1}\}$ are the basis
functions of a finite element defined using the functionals
$\{\modifiedfunctional_0,\dots,\modifiedfunctional_{n-1}\}$. Therefore,
by setting
$
  \mat{T}_{ij}
  =
  \left(\mat{B}^E_{\reflection}\right)^\alpha\left(\mat{B}^E_{\rotation}\right)^\beta
$
in \cref{eq:blocks}, we obtain a matrix $\mat{T}$ that corrects for
orientation differences on every sub-entity of $\refel$.

It follows from \cref{assume-me} that the base transformation matrices
for each sub-entity of the same type are the same. Therefore, once we
have computed the base transformations for one sub-entity of each type,
we can combine these to compute each block of the DOF transformation
matrix $\mat{T}$ for each cell. What remains is the determination of
$\alpha$ and $\beta$ for a given cell, which is covered in
\cref{sec:full_T}.

\section{Degree-of-freedom transformation algorithms}
\label{sec:algorithm}

In this section, we present our algorithm for computing the base
transformation matrices $\mat{B}^E_{\rotation}$ and
$\mat{B}^E_{\reflection}$.

\subsection{Base transformation algorithm}
\label{sec:computing-them}

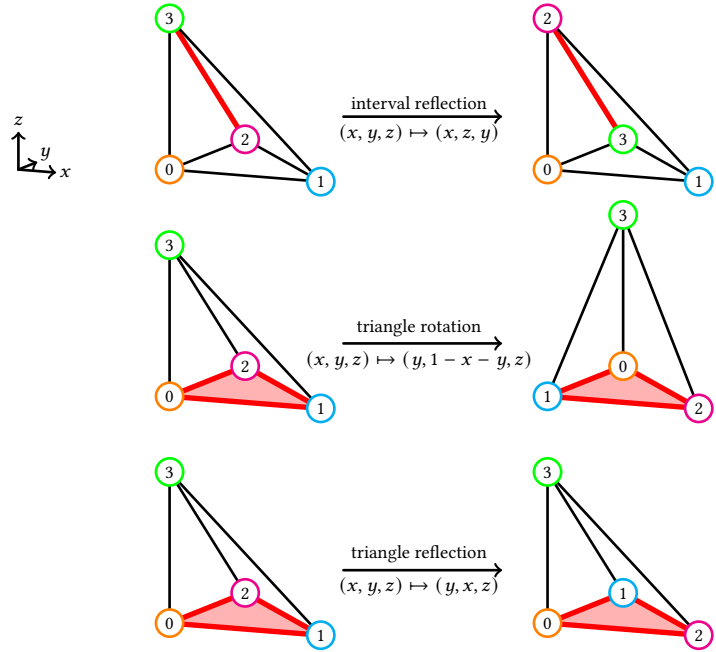
\begin{figure}
\begin{tikzpicture}[line width=1pt,line cap=round, line join=round]
\footnotesize
\newcommand{\pt}[3]{(#1+#2/2,-2*#1/25+#2/5+#3)}

\begin{scope}[shift={(-2,0)}]
\draw[->] \pt{0}{0}{0} -- \pt{0.5}{0}{0} node[anchor=west,inner sep=1pt] {$x$};
\draw[->] \pt{0}{0}{0} -- \pt{0}{0.5}{0} node[anchor=south west,inner sep=0pt] {$y$};
\draw[->] \pt{0}{0}{0} -- \pt{0}{0}{0.5} node[anchor=south,inner sep=2pt] {$z$};
\end{scope}

\begin{scope}[shift={(8.6,0)}]
\draw[white,->] \pt{0}{0}{0} -- \pt{0}{0}{0.5} node[anchor=south,inner sep=2pt] {$z$};
\end{scope}

\draw[line width=2pt,red] \pt{0}{0}{2} -- \pt{0}{2}{0};
\draw \pt{0}{0}{0} -- \pt{2}{0}{0} -- \pt{0}{2}{0} -- cycle;
\draw \pt{0}{0}{2} -- \pt{2}{0}{0};
\draw \pt{0}{0}{2} -- \pt{0}{0}{0};
\draw[orange,fill=white] \pt{0}{0}{0} circle (1.8mm) node[black] {0};
\draw[cyan,fill=white] \pt{2}{0}{0} circle (1.8mm) node[black] {1};
\draw[magenta,fill=white] \pt{0}{2}{0} circle (1.8mm) node[black] {2};
\draw[green,fill=white] \pt{0}{0}{2} circle (1.8mm) node[black] {3};
\draw[->] (2.3,0.7) -- (4.4,0.7);
\node[anchor=south] at (3.3,0.7) {interval reflection};
\node[anchor=north] at (3.3,0.7) {$(x,y,z)\mapsto(x,z,y)$};
\begin{scope}[shift={(5,0)}]
\newcommand{\ptm}[3]{\pt{#1}{#3}{#2}}
\draw[line width=2pt,red] \ptm{0}{0}{2} -- \ptm{0}{2}{0};
\draw \ptm{0}{0}{0} -- \ptm{2}{0}{0} -- \ptm{0}{2}{0} -- cycle;
\draw \ptm{0}{0}{2} -- \ptm{2}{0}{0};
\draw \ptm{0}{0}{2} -- \ptm{0}{0}{0};
\draw[orange,fill=white] \ptm{0}{0}{0} circle (1.8mm) node[black] {0};
\draw[cyan,fill=white] \ptm{2}{0}{0} circle (1.8mm) node[black] {1};
\draw[magenta,fill=white] \ptm{0}{2}{0} circle (1.8mm) node[black] {2};
\draw[green,fill=white] \ptm{0}{0}{2} circle (1.8mm) node[black] {3};
\end{scope}

\begin{scope}[shift={(0,-3)}]
\draw[line width=2pt,red,fill=red!30!white] \pt{0}{0}{0} -- \pt{2}{0}{0} -- \pt{0}{2}{0} -- cycle;
\draw \pt{0}{0}{2} -- \pt{0}{2}{0};
\draw \pt{0}{0}{2} -- \pt{2}{0}{0};
\draw \pt{0}{0}{2} -- \pt{0}{0}{0};
\draw[orange,fill=white] \pt{0}{0}{0} circle (1.8mm) node[black] {0};
\draw[cyan,fill=white] \pt{2}{0}{0} circle (1.8mm) node[black] {1};
\draw[magenta,fill=white] \pt{0}{2}{0} circle (1.8mm) node[black] {2};
\draw[green,fill=white] \pt{0}{0}{2} circle (1.8mm) node[black] {3};
\draw[->] (2.3,0.7) -- (4.4,0.7);
\node[anchor=south] at (3.3,0.7) {triangle rotation};
\node[anchor=north] at (3.3,0.7) {$(x,y,z)\mapsto(y, 1-x-y,z)$};
\begin{scope}[shift={(5,0)}]
\newcommand{\ptm}[3]{\pt{#2}{\numexpr2-#1-#2\relax}{#3}}
\draw[line width=2pt,red,fill=red!30!white] \ptm{0}{0}{0} -- \ptm{2}{0}{0} -- \ptm{0}{2}{0} -- cycle;
\draw \ptm{0}{0}{2} -- \ptm{0}{2}{0};
\draw \ptm{0}{0}{2} -- \ptm{2}{0}{0};
\draw \ptm{0}{0}{2} -- \ptm{0}{0}{0};
\draw[orange,fill=white] \ptm{0}{0}{0} circle (1.8mm) node[black] {0};
\draw[cyan,fill=white] \ptm{2}{0}{0} circle (1.8mm) node[black] {1};
\draw[magenta,fill=white] \ptm{0}{2}{0} circle (1.8mm) node[black] {2};
\draw[green,fill=white] \ptm{0}{0}{2} circle (1.8mm) node[black] {3};
\end{scope}
\end{scope}

\begin{scope}[shift={(0,-6)}]
\draw[line width=2pt,red,fill=red!30!white] \pt{0}{0}{0} -- \pt{2}{0}{0} -- \pt{0}{2}{0} -- cycle;
\draw \pt{0}{0}{2} -- \pt{0}{2}{0};
\draw \pt{0}{0}{2} -- \pt{2}{0}{0};
\draw \pt{0}{0}{2} -- \pt{0}{0}{0};
\draw[orange,fill=white] \pt{0}{0}{0} circle (1.8mm) node[black] {0};
\draw[cyan,fill=white] \pt{2}{0}{0} circle (1.8mm) node[black] {1};
\draw[magenta,fill=white] \pt{0}{2}{0} circle (1.8mm) node[black] {2};
\draw[green,fill=white] \pt{0}{0}{2} circle (1.8mm) node[black] {3};
\draw[->] (2.3,0.7) -- (4.4,0.7);
\node[anchor=south] at (3.3,0.7) {triangle reflection};
\node[anchor=north] at (3.3,0.7) {$(x,y,z)\mapsto(y,x,z)$};
\begin{scope}[shift={(5,0)}]
\newcommand{\ptm}[3]{\pt{#2}{#1}{#3}}
\draw[line width=2pt,red,fill=red!30!white] \ptm{0}{0}{0} -- \ptm{2}{0}{0} -- \ptm{0}{2}{0} -- cycle;
\draw \ptm{0}{0}{2} -- \ptm{0}{2}{0};
\draw \ptm{0}{0}{2} -- \ptm{2}{0}{0};
\draw \ptm{0}{0}{2} -- \ptm{0}{0}{0};
\draw[orange,fill=white] \ptm{0}{0}{0} circle (1.8mm) node[black] {0};
\draw[cyan,fill=white] \ptm{2}{0}{0} circle (1.8mm) node[black] {1};
\draw[magenta,fill=white] \ptm{0}{2}{0} circle (1.8mm) node[black] {2};
\draw[green,fill=white] \ptm{0}{0}{2} circle (1.8mm) node[black] {3};
\end{scope}
\end{scope}

\end{tikzpicture}
\caption{Three functions
that transform the sub-entities of a tetrahedron. In each case, the
sub-entity that the function is transforming is shown in red. In the
entity numbering used by Basix (see
\cref{fig:tet-reference}), the edge in the top plot is numbered 0 and
the face in the other two plots is numbered 3. In the entity numbering used
by Symfem, the edge and face are numbered 5 and 0 (respectively).
Note that each function shown here is equivalent to applying one of the generators in
\cref{table:cells} to the highlighted sub-entity.}
\label{fig:tetrahedron-maps}
\end{figure}

The algorithm for computing the base transformation matrices for the
sub-entity $E$ is presented in \cref{algorithm:1}. The inputs to the
algorithm are the reference basis functions $\hat{\phi}_a, \dots,
\hat{\phi}_b$ and functionals $\hat{l}_a, \dots, \hat{l}_b$ associated
with $E$, and the geometric maps $\mapsymbol^E_{\rotation}$ and
$\mapsymbol^E_{\reflection}$ as defined in
\cref{sec:sub-entity-transforming}. These geometric maps depend only on
the reference cell type (and importantly do not depend on the finite
element type), and so can be hard-coded for each reference cell shape.
For each transformation $\mapsymbol$ of the sub-entity $E$, the
algorithm applies the corresponding push-forward
$\mathcal{F}_\mapsymbol$ to the basis functions $\hat{\phi}_a, \dots,
\hat{\phi}_b$ on the reference cell $\refel$ to obtain the functions
$\phi_a, \dots, \phi_b$ on the image $\mapsymbol(\refel)$
(\cref{algorithm:1}, lines 3--4). The cell $\mapsymbol(\refel)$ does not
necessarily coincide with the reference cell $\refel$ (for example, see
the second function in \cref{fig:tetrahedron-maps}), although the image
of the sub-entity $\mapsymbol(E)$ will coincide with $E$.
The effects of the three transformations for the tetrahedron are shown
in \cref{fig:tetrahedron-maps}.
The algorithm next computes the base transformation
$\mat{B}^E_\mapsymbol$ corresponding to $\mapsymbol$
(\cref{algorithm:1}, lines 5--10): the value in the $i$th row and $j$th
column of the matrix is the value obtained when the functional
$\hat{l}_{a+j}$ is applied to the pushed forward basis function
$\phi_{a+i}$. Due to \cref{assume-me2}, the evaluation of this
functional only requires the value of the function restricted to $E$,
and so the fact that $\mapsymbol(\refel)$ and $\refel$ may not coincide
does not lead to evaluation of the function $\phi_{a+i}$ outside of its
cell.
\begin{algorithm}
\newcommand{\VARIABLE}[1]{\texttt{#1}}
\newcommand{\LET}{\textbf{let }}
\newcommand{\INPUT}{\textbf{input }}
\caption{Computing the base transformation matrices for the sub-entity $E$.}
\label{algorithm:1}
\begin{algorithmic}[1]
\STATE \INPUT $\{\hat{\phi}_a, \dots, \hat{\phi}_b\}$, $\{\hat{l}_a, \dots, \hat{l}_b\}$, $\mathbb{H}=\begin{cases}
\left\{\mapsymbol^E_{\rotation},\mapsymbol^E_{\reflection}\right\}&\text{$E$ is 2-dimensional}\\[4pt]
\left\{\mapsymbol^E_{\reflection}\right\}&\text{$E$ is 1-dimensional}
\end{cases}$
\FOR{$\mapsymbol\in\mathbb{H}$}
\STATE \LET $\mathcal{F}_{\mapsymbol}\text{ be the push-forward associated with the geometric map }\mapsymbol$
\STATE \LET $\{\phi_a, \dots, \phi_b\} = \{\mathcal{F}_{\mapsymbol}(\hat{\phi}_a), \dots, \mathcal{F}_{\mapsymbol}(\hat{\phi}_b)\}$
\STATE \LET $\mat{B}^E_{\mapsymbol}\in\mathbb{R}^{n_E\times n_E}$
\FOR{$i\in\{0,\dots,n_E-1\}$}
\FOR{$j\in\{0,\dots,n_E-1\}$}
\STATE $\left[\mat{B}^E_{\mapsymbol}\right]_{ij}\gets \hat{l}_{a+j}(\phi_{a+i})$
\ENDFOR
\ENDFOR
\ENDFOR
\RETURN $\left\{\mat{B}^E_H\,\middle|\,H\in\mathbb{H}\right\}$
\end{algorithmic}
\end{algorithm}
In \cref{algorithm:1}, base transformations are computed for one
sub-entity type only. It can easily be adapted with additional loop over
sub-entity types.

As an example, three functions
$\mapsymbol^{\interval}_{\reflection}$,
$\mapsymbol^{\rtriangle}_{\rotation}$ and
$\mapsymbol^{\rtriangle}_{\reflection}$
that could make up $\mathbb{H}$ for an edge and a face of a tetrahedron are shown in
\cref{fig:tetrahedron-maps}.

\subsection{Computing the base transformations with functionals defined using quadrature}

In practice, functionals are evaluated using quadrature, as described in
\cref{sec:basix-functionals}. An updated version of \cref{algorithm:1} that uses
functionals defined using quadrature is presented in
\cref{algorithm:2}, where $\{\hat{\vec{p}}_0, \dots,
\hat{\vec{p}}_{n_p-1}\}$ and $\mat{M}$ are the points and weights,
respectively, used to evaluate the functionals $\hat{l}_a, \dots,
\hat{l}_b$ that are associated with~$E$. This form of the algorithm is
suited to high-performance implementations, and a fast implementation is
included in Basix \cite{basix}.
Further details of this algorithms's implementation in Basix can be found in this
paper's software component \cite{supplement}.

\begin{algorithm}
\newcommand{\VARIABLE}[1]{\texttt{#1}}
\newcommand{\LET}{\textbf{let }}
\newcommand{\INPUT}{\textbf{input }}
\caption{Computing the base transformation matrices for the sub-entity
$E$ using numerical functionals.}
\label{algorithm:2}
\begin{algorithmic}[1]
\STATE \INPUT $\{\hat{\phi}_a, \dots, \hat{\phi}_b\}$, $\{\hat{\vec{p}}_0,\dots,\hat{\vec{p}}_{n_p-1}\}$, $\mat{M}$, $\mathbb{H}=\begin{cases}
\left\{\mapsymbol^E_{\rotation},\mapsymbol^E_{\reflection}\right\}&\text{$E$ is 2-dimensional}\\[4pt]
\left\{\mapsymbol^E_{\reflection}\right\}&\text{$E$ is 1-dimensional}
\end{cases}$
\FOR{$\mapsymbol\in\mathbb{H}$}
\STATE \LET $\mathcal{F}_\mapsymbol = \text{push-forward function associated with the map }\mapsymbol$
\STATE \LET $\{\vec{p}_0, \dots, \vec{p}_{n_p-1}\} = \{\mapsymbol^{-1}\hat{\vec{p}}_0, \dots, \mapsymbol^{-1}\hat{\vec{p}}_{n_p-1}\}$
\STATE \LET $\mat{B}^E_\mapsymbol \in \mathbb{R}^{n_E \times n_E}$
\FOR{$i \in \{0, \dots, n_E- 1 \}$}
\STATE \LET $\{v_0, \dots, v_{n_p-1}\} = \{\mathcal{F}_\mapsymbol (\phi_{a+i}(\vec{p}_0)), \dots, \mathcal{F}_\mapsymbol (\phi_{a+i}(\vec{p}_{n_p-1}))\}$
\FOR{$j\in\{0,\dots,n_E-1\}$}
\STATE $\left[\mat{B}^E_\mapsymbol\right]_{ij}\gets\displaystyle\sum_{l=0}^{s-1}\sum_{m=0}^{n_p-1} \mat{M}_{jlm}[v_m]_i$
\ENDFOR
\ENDFOR
\ENDFOR
\RETURN $\left\{\mat{B}^E_H\,\middle|\,H\in\mathbb{H}\right\}$
\end{algorithmic}
\end{algorithm}

\subsection{Computing the inverse and transpose transformations}

In \cref{sec:dof-transformations-matrix}, it was introduced that
$\vec{\phi} = \mat{T} \hat{\vec{\phi}}_g$ and shown that $\hat{\vec{c}} =
\mat{T}^{T} \vec{c}$, where $\hat{\vec{c}}$ is a vector of
degree-of-freedom values on a cell following the reference ordering and
$\vec{c}$ is the degrees-of-freedom following a globally consistent
ordering. To compute the inverse operations, we also requires the
inverse, transpose, and inverse transpose of~$\mat{T}$.

Computing the transpose of the base transformation matrices is trivial;
and once we have computed the inverse, computing the inverse transpose
is also trivial. We can use properties of the base transformation
matrices to compute the inverse base transformations without explicitly
inverting a matrix. If $\mat{B}^E_{\reflection}$ is the base
transformation matrix associated with reversing an edge, then we know
that
$$
\left(\mat{B}^{\interval}_{\reflection}\right)^2 = \mat{I}
$$
as reversing the edge twice is the same as doing nothing. Similarly, if
$\mat{B}^{\twodimentity}_{\rotation}$ and
$\mat{B}^{\twodimentity}_{\reflection}$ are the base transformation
matrices associated with rotating and reflecting a face with $n$ sides,
the we know that
\begin{align*}
  \left(\mat{B}^{\twodimentity}_{\rotation}\right)^n&=\mat{I},
  \\
  \left(\mat{B}^{\twodimentity}_{\reflection}\right)^2&=\mat{I},
\end{align*}
as reflecting the face twice or rotating $n$ times will return to the
original orientation. It follows that
\begin{align*}
\left(\mat{B}^{\interval}_{\reflection}\right)^{-1}
  &=\mat{B}^{\interval}_{\reflection},
\\
\left(\mat{B}^{\twodimentity}_{\rotation}\right)^{-1}
  &=\left(\mat{B}^{\twodimentity}_{\rotation}\right)^{n-1},
\\
\left(\mat{B}^{\twodimentity}_{\reflection}\right)^{-1}
  &=\mat{B}^{\twodimentity}_{\reflection},
\end{align*}
and so the inverse of each base transformation matrix is either equal to
the base transformation, or can be computed via $n-1$ matrix--matrix
multiplications.

\subsection{Computing the full DOF transformation matrix $\mat{T}$}
\label{sec:full_T}

Once the base transformation matrices have been computed, the full DOF
transformation matrix $\mat{T}$ can be computed by multiplying the base
transformations together to get each block. In practice, however, we
apply the base transformation matrices directly to the appropriate
subset of the basis functions for the cell rather than computing
explicit matrix--matrix products. In this section, we summarise how we
can determine the appropriate combination of matrices to multiply for
each sub-entity; this is covered in greater detail in \cite{2022-dofs}.
In this section, our discussion is based on using a low-to-high
orientation of sub-entities. Our method could be adapted if a different
orientation convention was used.

For each cell in a mesh, each vertex will have a local and a global
index. The global index is the vertex's common index across the full
mesh; the local index is the reference cell vertex that is mapped to the
`physical' vertex by the cell geometry map.

Consider an edge. If the vertex with the lowest global index does not
also have the lowest local index, the orientation of the edge is not
consistent with a low-to-high orientation.
In this case, the base transformation
$\mat{B}^{\interval}_{\reflection}$ should be applied to the basis
functions associated with the edge, or equivalently the block of
$\mat{T}$ for the edge is equal to $\mat{B}^{\interval}_{\reflection}$.
If the same vertex has both the lowest local and global indices, then
this block of $\mat{T}$ is the identity.

Consider a face. We define the vertex with the lowest global index to be
the `global origin' of the sub-entity, then look at the two neighbours
of this origin: a `global rotation' in the direction of the neighbour
with the lower global index is taken to be positive. We define the
vertex with the lowest local index to be the `local origin' of the
sub-entity, then look at the two neighbours of this origin: a `local
rotation' in the direction of the neighbour with the lower local index
is taken to be positive. If the global and local origins and/or rotation
directions do not agree, we determine the values of $\alpha$ and $\beta$
such that applying the transformation
$\left(G^{\twodimentity}_{\rotation}\right)^\beta \circ
\left(G^{\twodimentity}_{\reflection}\right)^\alpha$ to the
locally-numbered sub-entity gives a local orientation that matches the
global orientation. The block of the transformation matrix $\mat{T}$ for
this face will then be equal to
$\left(\mat{B}^{\twodimentity}_{\reflection}\right)^\alpha
\left(\mat{B}^{\twodimentity}_{\rotation}\right)^\beta$.

\subsection{Examples}

\begin{example}[Lagrange degree 3 on a quadrilateral]

The definition of this element can be found in
\cref{example-def-lagrange}. As a quadrilateral is a two-dimensional
cell, we only need to compute the base transformations for one of the
edges of the cell. There is a single base transformation
$\mat{B}^{\interval}_{\reflection}$ representing the effect of
reflecting a cell edge. To compute this base transformation, we use the
map $G^{\interval}_{\reflection} : (x,y) \mapsto (1-x, y)$. This map
will reverse edge 0 of the reference cell (using the numbering of the
reference sub-entities as shown in \cref{fig:quad-reference}). The basis
functions of the finite element that are associated with edge 0 are
\begin{align*}
\hat{\phi}_{4}(x,y) &= \tfrac{9}{4} x(1 - x)(2-3x) (1-y)(1-3y)(2-3y),
\\
\hat{\phi}_{5}(x,y) &= -\tfrac{9}{4} x(1 - x)(1-3x) (1-y)(1-3y)(2-3y).
\end{align*}
Applying the identity push-forward to these functions, we see that
\begin{align*}
\left[\mathcal{F}^\textup{id}(\hat{\phi}_{4})\right](x,y)
  &= \tfrac94 (1-x)x(3x-1) (1-y)(1-3y)(2-3y),
\\
\left[\mathcal{F}^\textup{id}(\hat{\phi}_{5})\right](x,y)
  &= -\tfrac94 (1-x)x(3x-2) (1-y)(1-3y)(2-3y).
\end{align*}
The functionals $\hat{l}_4$ and $\hat{l}_5$ associated with edge 0 are
point evaluations at the points $(\tfrac13, 0)$ and $(\tfrac23, 0)$.
Applying these to the mapped functions gives
\begin{align*}
\hat{l}_4\left(\mathcal{F}^\textup{id}(\hat{\phi}_{4})\right)
    &= \tfrac94\times\tfrac23\times\tfrac13\times0\times1\times1\times2=0,
    &\hat{l}_5\left(\mathcal{F}^\textup{id}(\hat{\phi}_{4})\right)
    &= \tfrac94\times\tfrac13\times\tfrac23\times1\times1\times1\times2=1,
\\
\hat{l}_4\left(\mathcal{F}^\textup{id}(\hat{\phi}_{5})\right)
    &= -\tfrac94\times\tfrac23\times\tfrac13\times-1\times1\times1\times2=1
    &\hat{l}_5\left(\mathcal{F}^\textup{id}(\hat{\phi}_{5})\right)
    &= -\tfrac94\times\tfrac13\times\tfrac23\times0\times1\times1\times2=0,
\end{align*}
and so
$$
  \mat{B}^{\interval}_{\reflection}
  =
  \begin{bmatrix}
  0&1\\1&0
  \end{bmatrix}.
$$
This base transformation matrix swaps the basis functions $\hat{\phi}_4$
and $\hat{\phi}_5$, which is what we would expect to happen if we
reversed the edge.
\end{example}

\begin{example}[N\'ed\'elec degree 2 on a tetrahedron]
The definition of this element can be found in
\cref{example-def-nedelec}. A tetrahedron is a three-dimensional cell,
so we must compute the base transformations for both an edge and a face
of the cell. We use the same maps as shown in
\cref{fig:tetrahedron-maps}, i.e.
\begin{align*}
G_{\reflection}^{\interval}&:(x,y,z)\mapsto(x,z,y),\\
G_{\rotation}^{\rtriangle}&:(x,y,z)\mapsto(y,1-x-y,z),\\
G_{\reflection}^{\rtriangle}&:(x,y,z)\mapsto(y,x,z).
\end{align*}
These maps will lead to the base transformations
$\mat{B}_{\reflection}^{\interval}$, $\mat{B}_{\rotation}^{\rtriangle}$,
and $\mat{B}_{\reflection}^{\rtriangle}$, respectively.

The map $G_{\reflection}^{\interval}$ reverses edge 0 of the reference
cell (using the numbering of the reference sub-entities as shown in
\cref{fig:tet-reference}). The basis functions associated with edge 0
are
\begin{align*}
  \hat{\phi}_0(x,y,z)
    &= \begin{bmatrix}
      0 \\ 2(1-4y)z \\ 4y(2y-1)
  \end{bmatrix},
  &\hat{\phi}_1(x,y,z)
  &=\begin{bmatrix}
    0\\4z(1-2z)\\2y(4z-1)
  \end{bmatrix}.
\end{align*}
Applying the covariant Piola push-forward map to these functions, we see
that
\begin{align*}
\left[\mathcal{F}_{\mapsymbol^{\scaleddown{\interval}}_{\reflection}}^\textup{curl}(\hat{\phi}_0)\right](x,y,z)
  &=\begin{bmatrix}
  1&0&0 \\
  0&0&1 \\
  0&1&0
  \end{bmatrix}
  \begin{bmatrix}
    0\\2(1-4z)y\\4z(2z-1)
  \end{bmatrix}
  =
  \begin{bmatrix}
    0\\4z(2z-1)\\2(1-4z)y
  \end{bmatrix},
  \\
\left[\mathcal{F}_{\mapsymbol^{\scaleddown{\interval}}_{\reflection}}^\textup{curl}(\hat{\phi}_1)\right](x,y,z)
  &=\begin{bmatrix}
    1&0&0 \\ 0&0&1 \\ 0&1&0
  \end{bmatrix}
  \begin{bmatrix}
    0 \\ 4y(1-2y) \\ 2z(4y-1)
  \end{bmatrix}
  =\begin{bmatrix}
    0 \\ 2z(4y-1) \\ 4y(1-2y)
  \end{bmatrix}.
\end{align*}
The functionals $\hat{l}_0$ and $\hat{l}_1$ are associated with edge 0
and are defined in \cref{example-def-nedelec}. Applying these to the
mapped functions gives
\begin{align*}
\hat{l}_0\left(\mathcal{F}_{\mapsymbol^{\scaleddown{\interval}}_{\reflection}}^\textup{curl}(\hat{\phi}_0)\right)
  &= \int_0^1\begin{bmatrix}0\\4t(2t-1)\\2(1-4t)(1-t)\end{bmatrix} \cdot \begin{bmatrix}0
    \\
    t-1
    \\
    1-t
  \end{bmatrix} \, \mathrm{d}t = 0,
  &\hat{l}_1\left(\mathcal{F}_{\mapsymbol^{\scaleddown{\interval}}_{\reflection}}^\textup{curl}(\hat{\phi}_0)\right)
  &= \int_0^1\begin{bmatrix}0\\4t(2t-1)\\2(1-4t)(1-t)\end{bmatrix}\cdot\begin{bmatrix}0
    \\
    -t
    \\
    t
  \end{bmatrix}\,\mathrm{d}t=-1,
  \\
\hat{l}_0\left(\mathcal{F}_{\mapsymbol^{\scaleddown{\interval}}_{\reflection}}^\textup{curl}(\hat{\phi}_1)\right)
  &=
  \int_0^1\begin{bmatrix}0
    \\
    2t(3-4t)
    \\
    4(1-t)(2t-1)
  \end{bmatrix}
  \cdot
  \begin{bmatrix}
    0
    \\
    t-1
    \\
    1-t
  \end{bmatrix} \, \mathrm{d}t = -1,
&\hat{l}_1\left(\mathcal{F}_{\mapsymbol^{\scaleddown{\interval}}_{\reflection}}^\textup{curl}(\hat{\phi}_1)\right)
  &=
  \int_0^1\begin{bmatrix}
    0
    \\
    2t(3-4t)
    \\
    4(1-t)(2t-1)
  \end{bmatrix}
  \cdot
  \begin{bmatrix}
    0
    \\
    -t
    \\
    t
  \end{bmatrix}\,\mathrm{d}t=0,
\end{align*}
and so
$$
\mat{B}^{\interval}_{\reflection}
= \begin{bmatrix}
  0&-1\\-1&0
\end{bmatrix}.
$$

The maps $G_{\rotation}^{\rtriangle}$ and $G_{\reflection}^{\rtriangle}$
rotate and reflect face 3 of the reference cell. The basis functions
associated with face 3 are
\begin{align*}
  \hat{\phi}_{18}(x,y,z) &= \begin{bmatrix}8y(2-x-2y-2z)\\8x(x+2y+z-1)\\8xy\end{bmatrix},
  &\hat{\phi}_{19}(x,y,z) &= \begin{bmatrix}8y(2x+y+z-1)\\8x(2-2x-y-2z)\\8xy\end{bmatrix}.
\end{align*}
Applying the covariant Piola push-forward map to these functions, we see
that
\begin{align*}
\left[\mathcal{F}_{\mapsymbol^{\scaleddown{\rtriangle}}_{\rotation}}^\textup{curl}(\hat{\phi}_{18})\right](x,y,z)
  &= \begin{bmatrix}-1&1&0\\-1&0&0\\0&0&1\end{bmatrix}\begin{bmatrix}8x(1-x+y-2z)\\8(1-x-y)(x-y+z)\\8(1-x-y)x\end{bmatrix}
=\begin{bmatrix}8(y-z)(-1-x+y)\\-8x(1-x+y-2z)\\8(1-x-y)x\end{bmatrix},
\\
\left[\mathcal{F}_{\mapsymbol^{\scaleddown{\rtriangle}}_{\rotation}}^\textup{curl}(\hat{\phi}_{19})\right](x,y,z)&=\begin{bmatrix}-1&1&0\\-1&0&0\\0&0&1\end{bmatrix}\begin{bmatrix}8x(1-x-2y+z)\\8(1-x-y)(x+2y-2z)\\8(1-x-y)x\end{bmatrix}
=\begin{bmatrix}8(y-z)(2-x-2y)\\-8x(1-x-2y+z)\\8(1-x-y)x\end{bmatrix},
\\
\left[\mathcal{F}_{\mapsymbol^{\scaleddown{\rtriangle}}_{\reflection}}^\textup{curl}(\hat{\phi}_{18})\right](x,y,z)&=\begin{bmatrix}0&1&0\\1&0&0\\0&0&1\end{bmatrix}\begin{bmatrix}8x(2-y-2x-2z)\\8y(y+2x+z-1)\\8xy\end{bmatrix}
=\begin{bmatrix}8y(y+2x+z-1)\\8x(2-y-2x-2z)\\8xy\end{bmatrix},
\\
\left[\mathcal{F}_{\mapsymbol^{\scaleddown{\rtriangle}}_{\reflection}}^\textup{curl}(\hat{\phi}_{19})\right](x,y,z)&=\begin{bmatrix}0&1&0\\1&0&0\\0&0&1\end{bmatrix}\begin{bmatrix}8x(2y+x+z-1)\\8y(2-2y-x-2z)\\8xy\end{bmatrix}
=\begin{bmatrix}8y(2-2y-x-2z)\\8x(2y+x+z-1)\\8xy\end{bmatrix}.
\end{align*}
The functionals $\hat{l}_{18}$ and $\hat{l}_{19}$ are associated with
face 3 and are defined in \cref{example-def-nedelec}. Applying these to
the mapped functions gives
\begin{align*}
\hat{l}_{18}\left(\mathcal{F}_{\mapsymbol^{\scaleddown{\rtriangle}}_{\rotation}}^\textup{curl}(\hat{\phi}_{18})\right)&=\int_0^1\int_0^{1-t}8t(-1-s+t)\,\mathrm{d}s\,\mathrm{d}t=-1,&
\hat{l}_{19}\left(\mathcal{F}_{\mapsymbol^{\scaleddown{\rtriangle}}_{\rotation}}^\textup{curl}(\hat{\phi}_{18})\right)&=\int_0^1\int_0^{1-t}-8s(1-s+t)\,\mathrm{d}s\,\mathrm{d}t=-1,\\
\hat{l}_{18}\left(\mathcal{F}_{\mapsymbol^{\scaleddown{\rtriangle}}_{\rotation}}^\textup{curl}(\hat{\phi}_{19})\right)&=\int_0^1\int_0^{1-t}8t(2-s-2t)\,\mathrm{d}s\,\mathrm{d}t=1,&
\hat{l}_{19}\left(\mathcal{F}_{\mapsymbol^{\scaleddown{\rtriangle}}_{\rotation}}^\textup{curl}(\hat{\phi}_{19})\right)&=\int_0^1\int_0^{1-t}-8s(1-s-2t)\,\mathrm{d}s\,\mathrm{d}t=0,\\
\hat{l}_{18}\left(\mathcal{F}_{\mapsymbol^{\scaleddown{\rtriangle}}_{\reflection}}^\textup{curl}(\hat{\phi}_{18})\right)&=\int_0^1\int_0^{1-t}8t(t+2s-1)\,\mathrm{d}s\,\mathrm{d}t=0,&
\hat{l}_{19}\left(\mathcal{F}_{\mapsymbol^{\scaleddown{\rtriangle}}_{\reflection}}^\textup{curl}(\hat{\phi}_{18})\right)&=\int_0^1\int_0^{1-t}8s(2-t-2s)\,\mathrm{d}s\,\mathrm{d}t=1,\\
\hat{l}_{18}\left(\mathcal{F}_{\mapsymbol^{\scaleddown{\rtriangle}}_{\reflection}}^\textup{curl}(\hat{\phi}_{19})\right)&=\int_0^1\int_0^{1-t}8t(2-2t-s)\,\mathrm{d}s\,\mathrm{d}t=1,&
\hat{l}_{19}\left(\mathcal{F}_{\mapsymbol^{\scaleddown{\rtriangle}}_{\reflection}}^\textup{curl}(\hat{\phi}_{19})\right)&=\int_0^1\int_0^{1-t}8s(2t+s-1)\,\mathrm{d}s\,\mathrm{d}t=0,
\end{align*}
and so
\begin{align*}
  \mat{B}^{\rtriangle}_{\rotation}
  &=
  \begin{bmatrix}
    -1 & - 1
    \\
    1 & 0
  \end{bmatrix},
  &\mat{B}^{\rtriangle}_{\reflection}
  &=
  \begin{bmatrix}
  0 & 1
  \\
  1 &0
  \end{bmatrix}.
\end{align*}
\end{example}

\subsection{Comparison to current algorithms}

As our algorithm for computing DOF transformation matrices has been newly
developed to compute the matrices required by the method proposed in \cite{2022-dofs},
comparisions with current algorithms are difficult to make. Prior to implementing
this algorithm, the base transformation matrices were computed manually by
carefully considering the DOF functionals. This manual
computation was challenging and time-consuming for more advanced elements,
and so the development of this algorithm has greatly reduced the effort required
to implement new elements, as well as allowing us to support custom elements
defined by the user at runtime, without requiring the user to first manually
compute the base transformation matrices.

A similar algorithm for computing transformation matrices was independently
developed and proposed in \cite{collin-thesis}. The implementation of this
algorithm is not currently publicly available, but we expect that it would
compute the same matrices as our algorithm, up to the choice of rotation
direction of a two dimensional entity (see \cref{table:cells}).

\section{In-place application of transformations}
\label{sec:in-place}

It is not necessary to compute the full DOF transformation matrix
$\mat{T}$. As described in \cref{sec:split_by_entity}, a combination of
the base transformation matrices can be applied directly and without
forming the full transformation matrix for a cell. It is possible to
efficiently apply each of the base transformation matrices in-place. We
describe now how this is be done.

A special, but common, case is when all base transformation matrices are
permutations. In practice, the effect of a permutation base
transformations can be applied to the degree-of-freedom map rather than
to each local element matrices or vectors. For the more general case,
the transformation is applied to the local element matrices and vectors.

\subsection{Permutation}
\label{sec:in-place-permutation}

We present the in-place application of permutations, based on
\cite{perm-in-place}. A clear and simple description of this method is
given in \cite{medium_perms}.

A vector $\vec{p} = (p_i) \in \mathbb{N}^n$ represents a permutation of
$n$ items if the entries of $\vec{p}$ are the numbers 0 to $n-1$ (with
each number appearing exactly one). A vector $\vec{p}$ of this form is
equivalent to the permutation matrix $\mat{P} = (\tilde{p}_{ij})$ with
entries
$$
  \tilde{p}_{ij}
  =
  \begin{cases}
    1&p_i=j,\\
    0&\text{otherwise}.
  \end{cases}
$$
The precomputation step for a permutation vector $\vec{p}$ is given in
\cref{algorithm:prepare_permutation}. For each $i = 0, \dots, n - 1$, we
replace $p_i$ with $p_{p_i}$ while $p_i < i$. Once this is done, the
permutation can be applied by following the method in
\cref{algorithm:apply_permutation}: for each $i = 0, \dots, n-1$, we swap the
$i$th and $p_i$th items in the array we are permuting.

\begin{algorithm}
\newcommand{\VARIABLE}[1]{\texttt{#1}}
\newcommand{\LET}{\textbf{let }}
\newcommand{\INPUT}{\textbf{input }}
\caption{Precomputation step for a permutation, where $\mathcal{P}(n)$ is the set of permutations of the integers $\{0,\dots,n-1\}$.}
\label{algorithm:prepare_permutation}
\begin{algorithmic}[1]
\STATE \INPUT $\vec{p}\in\mathcal{P}(n)$
\FOR{$i\in\{0,\dots,n-1\}$}
\WHILE{$p_i<i$}
\STATE $p_i\gets p_{p_i}$
\ENDWHILE
\ENDFOR
\RETURN $\vec{p}$
\end{algorithmic}
\end{algorithm}

\begin{algorithm}
\newcommand{\VARIABLE}[1]{\texttt{#1}}
\newcommand{\LET}{\textbf{let }}
\newcommand{\INPUT}{\textbf{input }}
\caption{In-place application of a permutation, where $\vec{p}$ is the result of applying \cref{algorithm:prepare_permutation} to a permutation and $\vec{d}=(d_i)$ is a vector of
$n$ items in some set $\mathbb{F}$.}
\label{algorithm:apply_permutation}
\begin{algorithmic}[1]
\STATE \INPUT $\vec{p}\in\mathbb{N}^n$, $\vec{d}\in\mathbb{F}^n$
\FOR{$i\in\{0,\dots,n-1\}$}
\STATE $d_i\leftrightarrow d_{p_i}$
\ENDFOR
\RETURN $\vec{d}$
\end{algorithmic}
\end{algorithm}

\subsection{Matrix--vector product}

The presented in-place matrix--vector multiplication is based on
\cite{in-place-matvec} and utilises LU decomposition. Let $\mat{B} \in
\mathbb{R}^{n \times n}$ be a base transformation matrix, and let
$\mat{P}, \mat{L}, \mat{U} \in \mathbb{R}^{n \times n}$ be a permutation
matrix, lower triangular matrix with 1s on the diagonal, and upper
triangular matrix (respectively) such that $\mat{B}\transpose =
\mat{P}^{-1} \mat{L} \mat{U}$. We know that there exists a natural
number $n > 0$ such that $\mat{B}^n = \mat{I}$, and so
$\det(\mat{B}\transpose) = \det(\mat{B}) \not = 0$. From this it follows
that the diagonal entries of $\mat{U}$ are non-zero.

In the precomputation step for in-place matrix--vector multiplication,
we overwrite entries of the matrix $\mat{B}$ above the diagonal with the
entries of $\mat{L}\transpose$ and the entries on or below the diagonal
with the entries of $\mat{U}\transpose$, and we store
$\mat{P}\invtranspose=\mat{P}$ in the vector form given in the previous
section so that an in-place permutation can be computed later. The
functions \pyth{sgetrf} (single precision) and \pyth{dgetrf} (double
precision) in LAPACK \cite{lapack} will compute the LU factorisation
in place and return the permutation in the correct format; as LAPACK
expects matrices to be input in column-major format, the LU decomposition
of the transpose that we desire is naturally computed if we pass a row-major
matrix $\mat{B}$ into one of these functions.

Noting that $\mat{B} \vec{v} = \mat{U} \transpose \mat{L}
\transpose\mat{P} \vec{v}$, we see that the matrix--vector product can
be computed in place by first applying the permutation $\mat{P}$ to
$\vec{v}$. This can be done in-place as in
\cref{sec:in-place-permutation}. We can then multiply by
$\mat{L}\transpose$ by adding $\sum_{j=i+1}^{n-1} l_{ij} v_j$ to $v_i$
for each $i = 0, \dots, n-1$ (where $v_i$ and $v_j$ are the entries of
$\vec{v}$ after the permutation has been applied). As the terms that we
add to $v_i$ only depend on $v_j$ for $j > i$, changing the values in
$\vec{v}$ as we go does not affect the results as long as we start by
adding to $v_0$ and proceed in order. We can then multiply by
$\mat{U}\transpose$ in a similar way by multiplying $v_i$ by $u_{ii}$
then adding $\sum_{j=0}^{i-1} u_{ij} v_j$. This time, the terms that we
add to $v_i$ include $v_j$ for $j < i$, so we must start with $v_{n-1}$
and iterate backwards through $\vec{v}$. This process is shown in 
\cref{algorithm:apply_matrix}.

\begin{algorithm}
\newcommand{\VARIABLE}[1]{\texttt{#1}}
\newcommand{\LET}{\textbf{let }}
\newcommand{\INPUT}{\textbf{input }}
\caption{Computing the matrix--vector product $\mat{B}\vec{d}$ in-place. If $\mat{P}^{-1}\mat{L}\mat{U}$ is the LU decomposition of $\mat{B}\transpose$, then the inputs to this
algorithm are: the matrix $\tilde{\mat{B}}=(\tilde{b}_{ij})$ whose entries on and below the diagonal are $\mat{U}\transpose$ and entries above the diagonal are $\mat{L}\transpose$;
the vector $\vec{p}=(p_i)$ which is the result of applying \cref{algorithm:prepare_permutation} to the permutation $\mat{P}$; and the vector $\vec{d}=(d_i)$.
Note that the first loop in this algorithms is an application of \cref{algorithm:apply_permutation}}
\label{algorithm:apply_matrix}
\begin{algorithmic}[1]
\STATE \INPUT $\tilde{\mat{B}}\in\mathbb{R}^{n\times n}$, $\vec{p}\in\mathbb{N}^n$, $\vec{d}\in\mathbb{R}^n$
\FOR{$i\in\{0,\dots,n-1\}$}
\STATE $d_i\leftrightarrow d_{p_i}$
\ENDFOR
\FOR{$i\in\{0,\dots,n-1\}$}
\FOR{$j\in\{i+1,\dots,n-1\}$}
\STATE $d_i\gets d_i+\tilde{b}_{ij}d_j$
\ENDFOR
\ENDFOR
\FOR{$i\in\{n-1,\dots,0\}$}
\STATE $d_i\gets m_{ii}d_i$
\FOR{$j\in\{0,\dots,i-1\}$}
\STATE $d_i\gets d_i+\tilde{b}_{ij}d_j$
\ENDFOR
\ENDFOR
\RETURN $\vec{p}$
\end{algorithmic}
\end{algorithm}

\section{Concluding remarks}
\label{sec:outro}

The new algorithm developed in this paper automates the computation of
degree-of-freedom transformations that allow arbitrary degree finite
element basis functions to be computed on a common reference cell whilst
preserving the required continuity of the global finite element space.
It overcomes a long-standing challenge of how to simply support
arbitrary degree finite elements on general meshes. Support for high
degree finite element spaces is important on modern computer
architectures, where the use of higher degree basis functions allows a
greater fraction of the available hardware performance to be exploited.

Our algorithm supports the implementation of finite elements from the
basic mathematical definition of an element, allowing a huge range of
elements to be implemented with minimal code. An implementation of this
approach is available in the FEniCSx libraries. However, there are
elements in the literature that cannot yet be used within FEniCSx. In
particular, elements such as Hermite \cite{hermite}, Bell \cite{bell},
and Argyris \cite{Argyris} that include evaluations of derivatives in
their functionals can currently be defined in Basix, but are not
supported the other components of FEniCSx. We believe that the methods
proposed in this paper can be used with minimal modification for these
elements, however we have not verified this with experiments.

\begin{acks}
  Support for MWS and GNW from EPSRC (EP/S005072/1) and Rolls-Royce plc
  as part of the Strategic Partnership in Computational Science for
  Advanced Simulation and Modelling of Engineering Systems (ASiMoV), and
  support for GNW from UKRI (UKRI1300) and EPSRC (EP/W026635/1) is
  gratefully acknowledged.
\end{acks}
\bibliographystyle{ACM-Reference-Format}
\ifthenelse{\equal{\anonymous}{YES}}{
\bibliography{refs-anonymous}
}{ 
\bibliography{refs}


\begin{thebibliography}{28}


\ifx \showCODEN    \undefined \def \showCODEN     #1{\unskip}     \fi
\ifx \showDOI      \undefined \def \showDOI       #1{#1}\fi
\ifx \showISBNx    \undefined \def \showISBNx     #1{\unskip}     \fi
\ifx \showISBNxiii \undefined \def \showISBNxiii  #1{\unskip}     \fi
\ifx \showISSN     \undefined \def \showISSN      #1{\unskip}     \fi
\ifx \showLCCN     \undefined \def \showLCCN      #1{\unskip}     \fi
\ifx \shownote     \undefined \def \shownote      #1{#1}          \fi
\ifx \showarticletitle \undefined \def \showarticletitle #1{#1}   \fi
\ifx \showURL      \undefined \def \showURL       {\relax}        \fi
\providecommand\bibfield[2]{#2}
\providecommand\bibinfo[2]{#2}
\providecommand\natexlab[1]{#1}
\providecommand\showeprint[2][]{arXiv:#2}

\bibitem[Agelek et~al\mbox{.}(2017)]%
        {agelek:2017}
\bibfield{author}{\bibinfo{person}{Rainer Agelek}, \bibinfo{person}{Michael
  Anderson}, \bibinfo{person}{Wolfgang Bangerth}, {and}
  \bibinfo{person}{William~L. Barth}.} \bibinfo{year}{2017}\natexlab{}.
\newblock \showarticletitle{On orienting edges of unstructured two- and
  three-dimensional meshes}.
\newblock \bibinfo{journal}{\emph{ACM Trans. Math. Software}}
  \bibinfo{volume}{44}, \bibinfo{number}{1}, Article \bibinfo{articleno}{5}
  (\bibinfo{year}{2017}), \bibinfo{numpages}{22}~pages.
\newblock
\urldef\tempurl%
\url{https://doi.org/10.1145/3061708}
\showDOI{\tempurl}


\bibitem[Anderson et~al\mbox{.}(1999)]%
        {lapack}
\bibfield{author}{\bibinfo{person}{E. Anderson}, \bibinfo{person}{Z. Bai},
  \bibinfo{person}{C. Bischof}, \bibinfo{person}{S. Blackford},
  \bibinfo{person}{J. Demmel}, \bibinfo{person}{J. Dongarra},
  \bibinfo{person}{J. Du~Croz}, \bibinfo{person}{A. Greenbaum},
  \bibinfo{person}{S. Hammarling}, \bibinfo{person}{A. McKenney}, {and}
  \bibinfo{person}{D. Sorensen}.} \bibinfo{year}{1999}\natexlab{}.
\newblock \bibinfo{booktitle}{\emph{{LAPACK} Users' Guide}
  (\bibinfo{edition}{third} ed.)}.
\newblock \bibinfo{publisher}{Society for Industrial and Applied Mathematics}.
\newblock


\bibitem[{A}rgyris et~al\mbox{.}(1968)]%
        {Argyris}
\bibfield{author}{\bibinfo{person}{{J}ohn~{H}. {A}rgyris},
  \bibinfo{person}{{I}saac {F}ried}, {and} \bibinfo{person}{{D}ieter~{W}.
  {S}charpf}.} \bibinfo{year}{1968}\natexlab{}.
\newblock \showarticletitle{The {TUBA} Family of Plate Elements for the Matrix
  Displacement Method}.
\newblock \bibinfo{journal}{\emph{The Aeronautical Journal}}
  \bibinfo{volume}{72}, \bibinfo{number}{692} (\bibinfo{year}{1968}),
  \bibinfo{pages}{{701--709}}.
\newblock
\urldef\tempurl%
\url{https://doi.org/10.1017/S000192400008489X}
\showDOI{\tempurl}


\bibitem[{A}rnold and {A}wanou(2011)]%
        {serendipity}
\bibfield{author}{\bibinfo{person}{{D}ouglas~{N}. {A}rnold} {and}
  \bibinfo{person}{{G}erard {A}wanou}.} \bibinfo{year}{2011}\natexlab{}.
\newblock \showarticletitle{{T}he serendipity family of finite elements}.
\newblock \bibinfo{journal}{\emph{Foundations of Computational Mathematics}}
  \bibinfo{volume}{11}, \bibinfo{number}{3} (\bibinfo{year}{2011}),
  \bibinfo{pages}{{337--344}}.
\newblock
\urldef\tempurl%
\url{https://doi.org/10.1007/s10208-011-9087-3}
\showDOI{\tempurl}


\bibitem[Arnold et~al\mbox{.}(2005)]%
        {mapping2}
\bibfield{author}{\bibinfo{person}{Douglas~N. Arnold}, \bibinfo{person}{Daniele
  Boffi}, {and} \bibinfo{person}{Richard~S. Falk}.}
  \bibinfo{year}{2005}\natexlab{}.
\newblock \showarticletitle{Quadrilateral $H({\rm div})$ Finite Elements}.
\newblock \bibinfo{journal}{\emph{SIAM J. Numer. Anal.}}  \bibinfo{volume}{42}
  (\bibinfo{year}{2005}), \bibinfo{pages}{2429--2451}.
\newblock
\urldef\tempurl%
\url{https://doi.org/10.1137/S0036142903431924}
\showDOI{\tempurl}


\bibitem[Baratta et~al\mbox{.}(2023)]%
        {fenicsx}
\bibfield{author}{\bibinfo{person}{Igor~A. Baratta}, \bibinfo{person}{Joseph~P.
  Dean}, \bibinfo{person}{J{\o}rgen~S. Dokken}, \bibinfo{person}{Michal
  Habera}, \bibinfo{person}{Jack~S. Hale}, \bibinfo{person}{Chris~N.
  Richardson}, \bibinfo{person}{Marie~E. Rognes}, \bibinfo{person}{Matthew~W.
  Scroggs}, \bibinfo{person}{Nathan Sime}, {and} \bibinfo{person}{Garth~N.
  Wells}.} \bibinfo{year}{2023}\natexlab{}.
\newblock \bibinfo{title}{{DOLFINx}: The next generation {FEniCS} problem
  solving environment}.  (\bibinfo{year}{2023}).
\newblock
\urldef\tempurl%
\url{https://doi.org/10.5281/zenodo.10447666}
\showDOI{\tempurl}
\newblock
\shownote{preprint}.


\bibitem[{B}ell(1969)]%
        {bell}
\bibfield{author}{\bibinfo{person}{{K}olbein {B}ell}.}
  \bibinfo{year}{1969}\natexlab{}.
\newblock \showarticletitle{{A} refined triangular plate bending finite
  element}.
\newblock \bibinfo{journal}{\emph{Internat. J. Numer. Methods Engrg.}}
  \bibinfo{volume}{1}, \bibinfo{number}{1} (\bibinfo{year}{1969}),
  \bibinfo{pages}{{101--122}}.
\newblock
\urldef\tempurl%
\url{https://doi.org/10.1002/nme.1620010108}
\showDOI{\tempurl}


\bibitem[{B}rezzi et~al\mbox{.}(1985)]%
        {bdm}
\bibfield{author}{\bibinfo{person}{{F}ranco {B}rezzi}, \bibinfo{person}{{J}im
  {D}ouglas}, {and} \bibinfo{person}{{L}.~{D}onatella {M}arini}.}
  \bibinfo{year}{1985}\natexlab{}.
\newblock \showarticletitle{{T}wo families of mixed finite elements for second
  order elliptic problems}.
\newblock \bibinfo{journal}{\emph{Numer. Math.}}  \bibinfo{volume}{47}
  (\bibinfo{year}{1985}), \bibinfo{pages}{{217--235}}.
\newblock
\urldef\tempurl%
\url{https://doi.org/10.1007/BF01389710}
\showDOI{\tempurl}


\bibitem[Ciarlet(1978)]%
        {Ciarlet:1978}
\bibfield{author}{\bibinfo{person}{Philippe~G. Ciarlet}.}
  \bibinfo{year}{1978}\natexlab{}.
\newblock \bibinfo{booktitle}{\emph{The Finite Element Method for Elliptic
  Problems}}.
\newblock \bibinfo{publisher}{North-Holland}, \bibinfo{address}{Amsterdam}.
\newblock


\bibitem[{C}iarlet and {R}aviart(1972)]%
        {hermite}
\bibfield{author}{\bibinfo{person}{{P}hilippe~{G}. {C}iarlet} {and}
  \bibinfo{person}{{P}ierre-{A}rnaud {R}aviart}.}
  \bibinfo{year}{1972}\natexlab{}.
\newblock \showarticletitle{{I}nterpolation theory over curved elements, with
  applications to finite element methods}.
\newblock \bibinfo{journal}{\emph{Computer Methods in Applied Mechanics and
  Engineering}} \bibinfo{volume}{1}, \bibinfo{number}{2}
  (\bibinfo{year}{1972}), \bibinfo{pages}{{217--249}}.
\newblock
\urldef\tempurl%
\url{https://doi.org/10.1016/0045-7825(72)90006-0}
\showDOI{\tempurl}


\bibitem[Collin(2025)]%
        {collin-thesis}
\bibfield{author}{\bibinfo{person}{T. Collin}.}
  \bibinfo{year}{2025}\natexlab{}.
\newblock \emph{\bibinfo{title}{Automated Finite Elements}}.
\newblock \bibinfo{thesistype}{Master's\ thesis}.
  \bibinfo{school}{Massachusetts Institute of Technology}.
\newblock
\urldef\tempurl%
\url{https://doi.org/1721.1/164593}
\showDOI{\tempurl}


\bibitem[Eisenstat(2014)]%
        {in-place-matvec}
\bibfield{author}{\bibinfo{person}{David Eisenstat}.}
  \bibinfo{year}{2014}\natexlab{}.
\newblock \bibinfo{title}{Is there an algorithm to multiply square matrices
  in-place?}
\newblock \bibinfo{howpublished}{Stack Overflow}.
\newblock
\urldef\tempurl%
\url{https://stackoverflow.com/a/25451717}
\showURL{%
\tempurl}
\newblock
\shownote{[Online; accessed 21-September-2026]}.


\bibitem[Ern and Jean-Luc(2021)]%
        {ern-guermond}
\bibfield{author}{\bibinfo{person}{Alexandre Ern} {and}
  \bibinfo{person}{Guermond Jean-Luc}.} \bibinfo{year}{2021}\natexlab{}.
\newblock \bibinfo{booktitle}{\emph{Finite Elements I: Approximation and
  Interpolation}}.
\newblock \bibinfo{publisher}{Springer}.
\newblock
\urldef\tempurl%
\url{https://doi.org/10.1007/978-3-030-56341-7}
\showDOI{\tempurl}


\bibitem[{Fich, Faith E. and Munro, J. Ian and Poblete, Patricio V.}(1995)]%
        {perm-in-place}
\bibfield{author}{\bibinfo{person}{{Fich, Faith E. and Munro, J. Ian and
  Poblete, Patricio V.}}} \bibinfo{year}{1995}\natexlab{}.
\newblock \showarticletitle{Permuting in place}.
\newblock \bibinfo{journal}{\emph{SIAM J. Comput.}} \bibinfo{volume}{24},
  \bibinfo{number}{2} (\bibinfo{year}{1995}), \bibinfo{pages}{{266--278}}.
\newblock
\urldef\tempurl%
\url{https://doi.org/10.1137/S0097539792238649}
\showDOI{\tempurl}


\bibitem[Kirby(2018)]%
        {mapping3}
\bibfield{author}{\bibinfo{person}{Robert~C. Kirby}.}
  \bibinfo{year}{2018}\natexlab{}.
\newblock \showarticletitle{A general approach to transforming finite
  elements}.
\newblock \bibinfo{journal}{\emph{The SMAI journal of computational
  mathematics}}  \bibinfo{volume}{4} (\bibinfo{year}{2018}),
  \bibinfo{pages}{197--224}.
\newblock
\urldef\tempurl%
\url{https://doi.org/10.5802/smai-jcm.33}
\showDOI{\tempurl}


\bibitem[Kirby et~al\mbox{.}(2012)]%
        {commonandunusual}
\bibfield{author}{\bibinfo{person}{Robert~C. Kirby}, \bibinfo{person}{Anders
  Logg}, \bibinfo{person}{Marie~E. Rognes}, {and} \bibinfo{person}{Andy~R.
  Terrel}.} \bibinfo{year}{2012}\natexlab{}.
\newblock \showarticletitle{Common and Unusual Finite Elements}.
\newblock In \bibinfo{booktitle}{\emph{Automated Solution of Differential
  Equations by the Finite Element Method}},
  \bibfield{editor}{\bibinfo{person}{Anders Logg}, \bibinfo{person}{Kent-Andre
  Mardal}, {and} \bibinfo{person}{Garth~N. Wells}} (Eds.).
  \bibinfo{series}{Lecture Notes in Computational Science and Engineering},
  Vol.~\bibinfo{volume}{84}. \bibinfo{publisher}{Springer},
  \bibinfo{address}{Heidelberg}, Chapter~3, \bibinfo{pages}{95--119}.
\newblock
\urldef\tempurl%
\url{https://doi.org/10.1007/978-3-642-23099-8_3}
\showDOI{\tempurl}


\bibitem[{N}\'ed\'elec(1980)]%
        {nedelec}
\bibfield{author}{\bibinfo{person}{{J}ean-{C}laude {N}\'ed\'elec}.}
  \bibinfo{year}{1980}\natexlab{}.
\newblock \showarticletitle{{M}ixed finite elements in \(\mathbb{{R}}^3\)}.
\newblock \bibinfo{journal}{\emph{Numer. Math.}} \bibinfo{volume}{35},
  \bibinfo{number}{3} (\bibinfo{year}{1980}), \bibinfo{pages}{{315--341}}.
\newblock
\urldef\tempurl%
\url{https://doi.org/10.1007/BF01396415}
\showDOI{\tempurl}


\bibitem[{N}\'ed\'elec(1986)]%
        {nedelec2}
\bibfield{author}{\bibinfo{person}{{J}ean-{C}laude {N}\'ed\'elec}.}
  \bibinfo{year}{1986}\natexlab{}.
\newblock \showarticletitle{{A} new family of mixed finite elements in
  \(\mathbb{{R}}^3\)}.
\newblock \bibinfo{journal}{\emph{Numer. Math.}} \bibinfo{volume}{50},
  \bibinfo{number}{1} (\bibinfo{year}{1986}), \bibinfo{pages}{{57--81}}.
\newblock
\urldef\tempurl%
\url{https://doi.org/10.1007/BF01389668}
\showDOI{\tempurl}


\bibitem[{R}aviart and {T}homas(1977)]%
        {rt}
\bibfield{author}{\bibinfo{person}{{P}ierre-{A}rnaud {R}aviart} {and}
  \bibinfo{person}{{J}ean-{M}arie {T}homas}.} \bibinfo{year}{1977}\natexlab{}.
\newblock \showarticletitle{{A} mixed finite element method for 2nd order
  elliptic problems}.
\newblock In \bibinfo{booktitle}{\emph{{M}athematical aspects of finite element
  methods}}, \bibfield{editor}{\bibinfo{person}{{I}lio {G}alligani} {and}
  \bibinfo{person}{{E}nrico {M}agenes}} (Eds.). Vol.~\bibinfo{volume}{606}.
  \bibinfo{pages}{{292--315}}.
\newblock


\bibitem[Rognes et~al\mbox{.}(2009)]%
        {rognes:2009}
\bibfield{author}{\bibinfo{person}{Marie~E. Rognes}, \bibinfo{person}{Robert~C.
  Kirby}, {and} \bibinfo{person}{Anders Logg}.}
  \bibinfo{year}{2009}\natexlab{}.
\newblock \showarticletitle{Efficient assembly of $H(\operatorname{div})$ and
  $H(\operatorname{curl})$ conforming finite elements}.
\newblock \bibinfo{journal}{\emph{SIAM Journal on Scientific Computing}}
  \bibinfo{volume}{31}, \bibinfo{number}{6} (\bibinfo{year}{2009}),
  \bibinfo{pages}{4130--4151}.
\newblock
\urldef\tempurl%
\url{https://doi.org/10.1137/08073901X}
\showDOI{\tempurl}


\bibitem[Scroggs(2021)]%
        {symfem}
\bibfield{author}{\bibinfo{person}{Matthew~W. Scroggs}.}
  \bibinfo{year}{2021}\natexlab{}.
\newblock \showarticletitle{Symfem: a symbolic finite element definition
  library}.
\newblock \bibinfo{journal}{\emph{Journal of Open Source Software}}
  \bibinfo{volume}{6}, \bibinfo{number}{64} (\bibinfo{year}{2021}),
  \bibinfo{pages}{{3556}}.
\newblock
\urldef\tempurl%
\url{https://doi.org/10.21105/joss.03556}
\showDOI{\tempurl}


\bibitem[Scroggs(2022)]%
        {basix-pr}
\bibfield{author}{\bibinfo{person}{Matthew~W. Scroggs}.}
  \bibinfo{year}{2022}\natexlab{}.
\newblock \bibinfo{title}{Compute entity transformations algorithmically for
  all elements rather than implementing them for each element (\#406)}.
\newblock \bibinfo{howpublished}{Github}.
\newblock
\urldef\tempurl%
\url{https://github.com/FEniCS/basix/pull/406}
\showURL{%
\tempurl}
\newblock
\shownote{[Online; accessed 21-September-2026]}.


\bibitem[Scroggs et~al\mbox{.}(2022a)]%
        {basix}
\bibfield{author}{\bibinfo{person}{Matthew~W. Scroggs},
  \bibinfo{person}{Igor~A. Baratta}, \bibinfo{person}{Chris~N. Richardson},
  {and} \bibinfo{person}{Garth~N. Wells}.} \bibinfo{year}{2022}\natexlab{a}.
\newblock \showarticletitle{Basix: a runtime finite element basis evaluation
  library}.
\newblock \bibinfo{journal}{\emph{Journal of Open Source Software}}
  \bibinfo{volume}{7}, \bibinfo{number}{73} (\bibinfo{year}{2022}),
  \bibinfo{pages}{{3982}}.
\newblock
\urldef\tempurl%
\url{https://doi.org/10.21105/joss.03982}
\showDOI{\tempurl}


\bibitem[Scroggs et~al\mbox{.}(2026a)]%
        {defelement-paper}
\bibfield{author}{\bibinfo{person}{Matthew~W. Scroggs},
  \bibinfo{person}{Pablo~D. Brubeck}, \bibinfo{person}{Joseph~P. Dean},
  \bibinfo{person}{J{\o}rgen~S. Dokken}, {and} \bibinfo{person}{India
  Marsden}.} \bibinfo{year}{2026}\natexlab{a}.
\newblock \showarticletitle{{DefElement:} an encyclopedia of finite element
  definitions}.
\newblock \bibinfo{journal}{\emph{Computational Science and Engineering}}
  \bibinfo{volume}{3} (\bibinfo{year}{2026}), \bibinfo{pages}{{2}}.
\newblock
\urldef\tempurl%
\url{https://doi.org/10.1007/s44207-026-00011-0}
\showDOI{\tempurl}


\bibitem[Scroggs et~al\mbox{.}(2026b)]%
        {defelement}
\bibfield{author}{\bibinfo{person}{Matthew~W. Scroggs},
  \bibinfo{person}{Pablo~D. Brubeck}, \bibinfo{person}{Joseph~P. Dean},
  \bibinfo{person}{J{\o}rgen~S. Dokken}, \bibinfo{person}{India Marsden},
  \bibinfo{person}{Nuno Nobre}, {et~al\mbox{.}}}
  \bibinfo{year}{{2020--2026}}\natexlab{b}.
\newblock \bibinfo{title}{{DefElement}: an encyclopedia of finite element
  definitions}.
\newblock \bibinfo{howpublished}{\url{https://defelement.org}}.
\newblock
\newblock
\shownote{[Online; accessed 21-September-2026]}.


\bibitem[Scroggs et~al\mbox{.}(2022b)]%
        {2022-dofs}
\bibfield{author}{\bibinfo{person}{Matthew~W. Scroggs},
  \bibinfo{person}{J{\o}rgen~S. Dokken}, \bibinfo{person}{Chris~N. Richardson},
  {and} \bibinfo{person}{Garth~N. Wells}.} \bibinfo{year}{2022}\natexlab{b}.
\newblock \showarticletitle{Construction of arbitrary order finite element
  degree-of-freedom maps on polygonal and polyhedral cell meshes}.
\newblock \bibinfo{journal}{\emph{ACM Trans. Math. Software}}
  \bibinfo{volume}{48}, \bibinfo{number}{2} (\bibinfo{year}{2022}),
  \bibinfo{pages}{{18:1--18:23}}.
\newblock
\urldef\tempurl%
\url{https://doi.org/10.1145/3524456}
\showDOI{\tempurl}


\bibitem[Scroggs and Wells(2026)]%
        {supplement}
\bibfield{author}{\bibinfo{person}{Matthew~W. Scroggs} {and}
  \bibinfo{person}{Garth~N. Wells}.} \bibinfo{year}{2026}\natexlab{}.
\newblock \bibinfo{title}{Software component for Algorithm XXXX: Computation of
  finite element degree-of-freedom transformation matrices}.
\newblock
\newblock
\urldef\tempurl%
\url{https://doi.org/10.5281/zenodo.15363224}
\showDOI{\tempurl}


\bibitem[user \href{https://medium.com/@kevingxyz}{@kevingxyz})(2019)]%
        {medium_perms}
\bibfield{author}{\bibinfo{person}{Kevin~(Medium user
  \href{https://medium.com/@kevingxyz}{@kevingxyz})}.}
  \bibinfo{year}{2019}\natexlab{}.
\newblock \bibinfo{title}{Permutation in place}.
\newblock
  \bibinfo{howpublished}{\url{https://medium.com/@kevingxyz/permutation-in-place-8528581a5553}}.
\newblock
\newblock
\shownote{[Online; accessed 21-September-2026]}.


\end{thebibliography}
}
\end{document}